\documentclass[12pt]{amsart}
\usepackage{amsfonts}
\usepackage{amssymb,amsmath}
\usepackage{amsthm}
\usepackage{epsfig}
\usepackage{color}

\newtheorem{theorem}{Theorem}[section]
\newtheorem{proposition}[theorem]{Proposition}
\newtheorem{lemma}[theorem]{Lemma}
\newtheorem{corollary}[theorem]{Corollary}

\newtheorem{definition}[theorem]{Definition}
\newtheorem*{thm:1.2-2}{Theorem \ref{thm:1.2-1}}
\newtheorem*{cor:1.3-2}{Corollary \ref{cor:1.3-1}} 
\newtheorem*{thm:1.4-2}{Theorem \ref{thm:1.4-1}}
  
\newtheorem*{thm:1.6-2}{Theorem \ref{thm:1.6-1}} 
\newtheorem*{thm:1.7-2}{Theorem \ref{thm:1.7-1}} 
\newtheorem*{thm:1.8-2}{Theorem \ref{thm:1.8-1}} 
\newtheorem*{thm:1.9-2}{Theorem \ref{theorem:shorteqns}} 
\newtheorem*{thm:1.10-2}{Theorem \ref{fieldK}} 
\newtheorem*{cor:1.11-2}{Corollary \ref{cor:1.11-1}} 
\newtheorem*{thm:1.12-2}{Theorem \ref{thm:1.12-1}} 
\newtheorem*{cor:1.13-2}{Corollary \ref{cor:totreal}}
\newtheorem*{thm:1.14-2}{Theorem \ref{theorem:squares}}
\newtheorem*{thm:1.15-2}{Theorem \ref{thm:1.15-1}}

\DeclareMathOperator{\Iso}{Iso}
\DeclareMathOperator{\HAG}{Aff^+(H_1(S,\QQ),\omega)}
\DeclareMathOperator{\Hom}{Hom}
\DeclareMathOperator{\hol}{hol}
\DeclareMathOperator{\Ann}{Ann}

\def\RR{{\mathbb R}}
\def\CC{{\mathbb C}}

\def\PP{{\mathbb P}}
\def\QQ{{\mathbb Q}}

\def\slr{{SL(2,\RR)}}

\newlength{\figboxwidth}
\begin{document}

\centerline{\large Algebraically Periodic Translation Surfaces} \vspace{10mm}
\medskip

\centerline{Kariane Calta and John Smillie}

\begin{section}{Introduction}
A translation surface is a compact surface with a geometric structure. Translation surfaces arise in complex analysis, the classification of surface diffeomorphisms and the study of polygonal billiards. This structure can be defined via complex analytic methods (where it is called an abelian differential). It can also be described as a pair of measured foliations but for our purposes the following description is more convenient.

A translation surface can be obtained by gluing a finite number of polygons $P_1, \ldots, P_n \subset \RR^2$ along parallel edges in such a way that the result is a compact surface.

In this paper we synthesize and further develop an algebraic theory of translation surfaces that builds upon the work of many authors. ( To name a very few: \cite{gj1}, \cite{calta}, \cite{kenyon-smillie}, \cite{arnoux1}, \cite{arnoux2}, \cite{ctm:billiards}). Many of the results in this paper are of a purely algebraic nature but they are all motivated by a desire to understand the geometry and dynamics of translation surfaces.

Given a translation surface $S$ and an $\alpha\in\slr$ we obtain a new translation surface $\alpha S$ by letting $\alpha$ act on the polygons $P_1, \ldots, P_n$ that comprise $S$. The {\em Veech group} of $S$ consists of those $\alpha\in\slr$ for which $\alpha S$ is equivalent to $S$ as a translation surface. If the Veech group of $S$ is a lattice in $\slr$, then we say that $S$ is a {\em lattice surface}.

 Algebraic tools play a role in several recent advances in the study of translation surfaces.
Kenyon and Smillie \cite{kenyon-smillie} introduced a translation surface invariant, the $J$-invariant, which takes values in $\RR^2 \wedge_{\QQ} \RR^2$. They used the $J$-invariant as a tool in the classification ofl lattice surfaces that arise from acute, rational triangular billiard tables. (\cite{kenyon-smillie},\cite{puchta})

Calta (\cite{calta}) and McMullen (\cite{ctm:billiards}) both use algebraic methods in their independent work on genus two translation structures.
In \cite{calta}, the $J$-invariant and projections of the $J$-invariant are used to develop algebraic criteria for the flow on a surface in a given direction to be periodic. A projection of the $J$-invariant in a direction $v$ is in fact the $SAF$ invariant (\cite{arnoux1}) associated to an interval exchange transformation induced by the flow in the direction $v$. These ideas played an important role in her discovery of infinitely many lattice surfaces in genus two.
 
McMullen (\cite{mc1}) discovered these same surfaces by complex analytic methods but he does use methods similar in spirit to the $J$ invariant in connection with other work on genus two translation surfaces.  In particular he  uses the notions of the complex flux and the Galois flux to locate genus two surfaces whose Veech groups are infinitely generated (\cite{ctm:billiards}). 

In addition to finding lattice surfaces Calta also uses the $J$-invariant to locate other genus two surfaces  that display interesting dynamical behavior: these are the completely periodic surfaces.  On a completely periodic surface, if $SAF_v$, the projection of $J$ in the direction $v$, vanishes then the geodesic flow in that direction is periodic. A direction in which each orbit of the geodesic flow is either closed or a saddle connection is called {\em completely periodic}. Calta calls those directions $v$ for which $SAF_v=0$ {\em algebraically periodic}. Thus, on a completely periodic surface, every algebraically periodic direction is completely periodic. For these surfaces, dynamical behavior can be deduced from algebraic information.

Our first result makes a connection between algebraic periodicity and number fields.

\begin{theorem} \label{alg_per} If $S$ has three algebraically periodic directions
then it has infinitely many.
If we choose coordinates so that the directions with slopes 0, 1 and $\infty$ are algebraically
periodic then there is a  field $K$ so that the collection
of algebraically periodic directions are exactly those with
slopes in $K \cup \lbrace \infty \rbrace$. Moreover $K$ is a number field with $\deg(K)\le genus(S)$.
\end{theorem}

If a surface $S$ has three, hence infinitely many, algebraically periodic directions we will say that $S$ is {\em algebraically periodic} and we will call the field $K$ described above the {\em periodic direction field} of $S$. 

We next show that algebraically periodic surfaces arise in many contexts.

Suppose that a translation surface $S$ has a pseudo-Anosov diffeomorphism $A$ and let $\lambda$ be its expansion constant. By replacing $A$ with $A^2$ if necessary we may assume that $\lambda>0$. Recall that the trace field of $S$ is the field $\QQ[\lambda + \lambda^{-1}]$. 
Gutkin-Judge and Kenyon-Smillie 
show that the trace field is the holonomy field of the surface.  In particular the the trace field depends only on $S$ and not  on the particular $A$ chosen (\cite{gj}, \cite{kenyon-smillie}
).  

We say that $S$ is {\em completely algebraically periodic} (not to be confused with completely periodic defined above) if every homological direction is algebraically periodic (see Definition~\ref{def:hd}). In Corollary~\ref{cor:capeq} we show that this property is equivalent to the equality of the periodic direction field and the holonomy field.

\begin{theorem} \label{thm:1.2-1} If a translation surface $S$ has an affine pseudo-Anosov automorphism then $S$ is algebraically periodic and its periodic
direction field is its trace field.
\end{theorem}

Thurston (\cite{thurston}) described a construction of translation surfaces starting with a pair of multi-curves (see also Veech \cite{veech}). Such surfaces contain pseudo-Anosov maps in their Veech groups so we have the following.

\begin{corollary} \label{cor:1.3-1} If $S$ is obtained from the Thurston-Veech construction, for example if $S$ is a lattice surface, then $S$ is completely algebraically periodic.
\end{corollary}

The lattice property for translation surfaces arising from rational triangles has been studied by Veech (\cite{veech}), Ward (\cite{ward}), Vorobets (\cite{vorobets}), Kenyon-Smillie (\cite{kenyon-smillie}), Puchta (\cite{puchta}) and Hooper (\cite{hooper}).

\begin{theorem} \label{thm:1.4-1} If $S$ is a flat surface obtained from a rational-angled triangle via the Zemlyakov-Katok construction then $S$ is completely algebraically periodic.
\end{theorem}

McMullen defined and used the {\em homological affine group} in (\cite{ctm:dynamics}) to describe those genus two surfaces whose $\slr$ orbits fail to be dense in their respective strata. This notion was first introduced by Gutkin and Judge in (\cite{gj}). We recall the definition here. 
The group $\HAG$ consists of
matrices $\alpha\in\slr$ for which there is
an automorphism $A$ of $H_1(S,\QQ)$ which preserves the
intersection numbers of homology classes and such that:

$$
\begin{array}{cccc}
 H_1(S,\QQ) & \xrightarrow{\hol} & \RR^2 \\
\downarrow{A} &  & \downarrow{\alpha} \\
 H_1(S,\QQ) & \xrightarrow{\hol} & \RR^2\end{array}
 $$
Here {\em $\hol$} is the holonomy homomorphism as defined in Section~\ref{sec:basics}.

\begin{theorem}  \label{thm:1.8-1} If $S$ is completely algebraically periodic, in particular if the Veech group of $S$ contains a hyperbolic element, then $$\HAG=SL(2,K).$$ 
\end{theorem}

In (\cite{calta}) Calta gives explicit formulae in terms of parameters of cylinder decompositions which allow one to recognize the completely periodic surfaces in genus two. Here we give formulae which allow one to recognize algebraically periodic surfaces of arbitrary genera and compute their periodic direction fields in terms of a triangulation of the surface. 

Consider a translation surface $S$ with a triangulation dividing the
surface into triangles $\Delta_j$. Say that $v_j$ and $w_j$ are two
successive sides of $\Delta_j$ which correspond to a counterclockwise orientation of $\Delta_j$.  Say that the coordinates of  $v_j$ are $\begin{bmatrix} a_j\\ c_j\end{bmatrix}$
and those of $w_j$ are $\begin{bmatrix} b_j\\ d_j\end{bmatrix}$. 

\begin{theorem}
\label{thm:1.6-1} Assume that $a_j, b_j, c_j,$ and $d_j$ lie in a number field $L$. Then the lines of slope 0, 1 and $\infty$ are algebraically periodic directions if and only if  for any pair of distinct field embeddings $\sigma, \tau$ of $L$ into
$\CC$, the following equations hold: 

\[ \sum_j \sigma(c_j) \tau(b_j)- \sigma(b_j) \tau(c_j)= \sum_j \sigma(d_j) \tau(a_j)- \sigma(a_j) \tau(d_j)
\]

\[
\sum_j \sigma(c_j) \tau(d_j)- \sigma(d_j) \tau(c_j)= \sum_j
\sigma(a_j)\tau(b_j)-\sigma(b_j)\tau(a_j)=0
\]
Furthermore if the above equations hold, then the periodic direction field is the field consisting of $\lambda\in L$  for which $\sigma(\lambda)=\tau(\lambda)$ whenever $\sum_j \sigma(c_j) \tau(b_j)- \sigma(b_j) \tau(c_j)\ne0$.
\end{theorem}

It is not known which fields can be obtained as trace fields of surfaces whose Veech groups contain pseudo-Anosov automorphisms.  We show that there are no restrictions on possible periodic direction fields.

\begin{theorem} \label{fieldK}
Every number field $K$ is the periodic direction
field for a completely algebraically periodic translation surface that arises from a right-angled billiard table. 
\end{theorem}

Given a quadratic subfield $K$ of the reals, McMullen
gives examples of L-shaped tables with explicit side lengths so that the corresponding translation surfaces have genus two and periodic direction field equal to $K$ (\cite{mc1}). We can view the examples in the previous result as being related to these.

Recall that the periodic direction field $K$ is a number field. Under certain assumptions, even more can be said about $K$. 
In the following results, we show that certain geometric properties of $S$ imply that $K$ is totally real.
If the Veech group of a translation surface contains
two noncommuting parabolic elements, in particular if the surface is a lattice
surface, Hubert and Lanneau show that the holonomy
field of the surface is totally real (\cite{hl}). We prove a variation of this result in Corollary~\ref{cor:totreal}.

\begin{corollary} \label{cor:totreal} If $S$ has a parabolic automorphism and a second
algebraically periodic direction then $S$ is algebraically periodic
and the periodic direction field is totally real.
\end{corollary}

If $S$ has a parabolic automorphism then $S$ has a decomposition into cylinders of rationally related moduli. By adjusting the heights and widths of the cylinders we see that $S$ is affinely equivalent to a translation surface $S'$ with a decomposition into cylinders of rational moduli. 
The surface $S'$ has a decomposition into squares. (In general this will not be an edge to edge decomposition). We are interested in the general question of which translation surfaces have decompositions into squares

The question of when a rectangle can be tiled by squares was answered by Dehn (\cite{dehn}). Kenyon (\cite{kenyon}) addresses the same question for the torus.

\begin{theorem}\label{theorem:squares} Suppose that $S$ is algebraically periodic and that $S$ has a decomposition into squares so that the directions of the sides are algebraically periodic directions. Then the periodic direction field $K$ is totally real. \end{theorem}
 
We emphasize that the squares need not have the same side lengths and that the decomposition need not be an edge to edge decomposition.

 The following result is a partial converse.

\begin{theorem} \label{thm:1.15-1}
Given a totally real number field $K$ there is a completely algebraically periodic translation surface with periodic direction field $K$ which has a decomposition into squares. \end{theorem}

We say that translation surfaces $S$ and $S'$ are {\em scissors congruent} if both surfaces can be assembled from the same collection of triangles in $\RR^2$ with perhaps different collections of gluing maps. As usual we require that the gluing maps identify edges with edges and are restrictions of translations. If we did not require that the gluing maps identify entire edges then we would get a much weaker notion of scissors congruence. For this weaker notion the only scissors congruence invariant is area. For the stronger notion of scissors congruence that we use here there are additional invariants.

\begin{theorem}
The property of being algebraically periodic with periodic direction field $K$ is a scissors congruence invariant.
\end{theorem}

A fundamental scissors congruence invariant is the $J$ invariant of a surface.
The $J$ invariant of a surface, like the homological affine group, captures information about the holonomy and the product pairing on cohomology of the surface $S$. Unlike the homological affine group the $J$ invariant is a scissors congruence invariant.
We define a group $\Iso(J(S))$ which consists of those $\alpha\in SL(2,\RR)$ which preserve the $J$ invariant of $S$. The Veech group of $S$ is contained in the homological affine which is contained in $\Iso(J(S))$.

The Veech group of a surface may be quite difficult to calculate. Although Calta and McMullen constructed explicit examples of Veech surfaces in genus two (the L-tables), the Veech groups of most of these surfaces remain unknown though Bainbridge has computed a large number of individual examples. By contrast, $\Iso(J)$ is relatively easy to compute:

\begin{theorem} \label{thm:1.5-1} If $S$ is algebraically periodic  with periodic direction field $K$ then $\Iso(J)$ is conjugate to $SL(2,K)$.
\end{theorem}

The isomorphism type of $\Iso(J)$ as an abstract group determines both whether $S$ is algebraically periodic and, if so, its periodic direction field.

\begin{theorem} \label{thm:1.7-1} If the group $\Iso(J(S))$ is isomorphic as a group to the group $SL(2,K)$ for some field $K$ then $S$ is algebraically periodic and $K$ is its periodic direction field.
\end{theorem}

Though $\Iso(J(S))$ and $\HAG$ are closely related, they are not necessarily the same.  In Corollary \ref{cor:notthesame} we give an example of a translation surface for which $\HAG\ne\Iso(J)$.

A priori the $J$ invariant of a surface $S$ depends on $S$. We show that when $S$ is algebraically periodic there are
cases in which the $J$ invariant is determined by the periodic direction field of $S$. 

The following is a corollary of Theorem~\ref{lemma:uniqueJ}. 

\begin{corollary} \label{cor:1.11-1} If $S$ is completely algebraically periodic then the $J$ invariant of $S$ is completely determined by the periodic direction field and the area of $S$.
\end{corollary}

This result is used in the construction of examples of translation surfaces with a given periodic direction field.

We say that $J \in \RR^2 \wedge_{\QQ} \RR^2$ is {\em  algebraically periodic} if $SAF_v(J)=0$ in at least three distinct directions. 

The following result gives an explicit formula for $J$ in terms of the periodic direction field.

\begin{theorem} \label{thm:1.12-1} If $\alpha_1\ldots\alpha_n$ and $\beta_1\ldots\beta_n$ are  dual bases for the trace form on a field $K$ then 
$$J=\sum_j \begin{bmatrix}\alpha_j\\0\end{bmatrix}\wedge \begin{bmatrix}0\\ \beta_j\end{bmatrix}$$ 
is algebraically periodic with periodic direction field $K$. Conversely if $J$ can be written as above and with entries in $K$ and $J$ is algebraically periodic with periodic direction field $K$ and area 1 then $\alpha_1\ldots\alpha_n$ and $\beta_1\ldots\beta_n$ are dual bases for the trace form on $K$.
\end{theorem}

The questions involving square tilings of surfaces are related to the question of when $J$ be written in the form
$$J=\sum_j \begin{bmatrix}\alpha_j\\0\end{bmatrix}\wedge \begin{bmatrix}0\\ \beta_j\end{bmatrix}$$ with $\alpha_j=\beta_j$.

We connect this question with to the signature of aa symmetric bilinear form similar to that considered in (\cite{kenyon}). We show that the signature of this form is also related to the number of real and complex embeddings of the field $K$.

\end{section}

\begin{section}{Basics}
\label{sec:basics}

The $J$-invariant of a surface $S$ plays a central role in this paper. We recall the definition of this invariant  given in (\cite{kenyon-smillie}).
The $J$-invariant takes values in $\RR^2\wedge_\QQ\RR^2$. Specifically, if $\Delta$ is a triangle with vertices $p_0$, $p_1$ and $p_2$ labeled in counterclockwise order, then
$$J(\Delta)=p_0 \wedge p_1 + p_1 \wedge p_2 +p_2\wedge p_0.$$
Note that this quantity does not depend on which vertex is labeled $p_0$.
On the other hand, if $v$ and $w$ are the coordinate vectors corresponding to two of the sides of $\Delta$ and if the ordered pair $v,w$ determines the correct orientation then
$$J(\Delta)=v\wedge w.$$ In particular, the value of $J$ remains unchanged
if the triangle is translated.
Recall that $S$ has a cellular or edge-to-edge decomposition into triangles if neighboring triangles share entire edges.
If a surface $S$ has a cellular decomposition into triangles $\Delta_1\ldots \Delta_n$ then
we define $J(S)$ by the following formula:
\begin{equation}\label{Jdef}
J(S)= \sum_{j=1}^n J(\Delta_j).
\end{equation}
The resulting element is independent of the cellular decomposition chosen.
This remains true if we add a finite collection of points which can
serve as additional vertices in the cellular decomposition.

\begin{proposition} The $J$ invariant of a surface $S$ is a scissors congruence invariant.
\end{proposition}

We can also define the $J$ invariant for translation surfaces with boundaries as long as the boundaries are finite unions of segments. We can do this by adding interior vertices to the collection of singular points, triangulating and summing the contribution of the triangles. The resulting element of $\RR^2\wedge_\QQ\RR^2$ does not depend on the triangulation chosen.  Note though that while adding vertices to the interior of the surface does not change the value of the $J$ invariant, adding vertices to the boundary of the surface can change the value.
We refer readers to (\cite{kenyon-smillie}) for proofs of these facts.

The $J$-invariant  is closely related to the Sah-Arnoux-Fathi $\operatorname{SAF}$  invariants of the individual
direction flows. We recall the $\operatorname{SAF}$-invariant
 of an interval exchange map.
This invariant takes its values in $\RR\wedge_\QQ\RR$.
If the interval exchange takes intervals of lengths $\ell_1,\ldots\ell_n$ and translates
them by $t_1,\ldots t_n$ then the value of the invariant is
\[
\sum_{j=1}^n \ell_j\wedge t_j.
\]

 If $v\in\RR^2$ is a non-zero vector then there is a flow in the
direction $v$ and an $\operatorname{SAF}$-invariant associated to this flow which
we denote by $\operatorname{SAF}_v$. 
Let $\pi:\RR^2\to\RR$ be the orthogonal  projection with $v$ in its kernel. Then $\pi$
induces a map from $\RR^2\wedge_\QQ \RR^2\to \RR\wedge_\QQ \RR$
and the image of $J(S)$ under this map is $\operatorname{SAF}_v(S)$.

The following is 
Theorem $21$ of (\cite{kenyon-smillie}):
\begin{theorem} \label{thm:ks21} If $S$ has a cylinder decomposition in the
direction $v$, then $\operatorname{SAF}_v(S)=0$. \end{theorem}

In \cite{calta}, a direction $v$ was said to be {\em completely periodic} if $S$ decomposed as a union of cylinders in the direction $v$. We are interested in
a related notion.

\begin{definition} A direction $\bar{v}$ on a surface $S$ is said to be
{\it algebraically periodic} if $\operatorname{SAF}_{\bar{v}}(S)=0$. \end{definition}

According to Theorem~\ref{thm:ks21}  any completely periodic direction is algebraically periodic but the converse is not true.

We recall several basic notions related to translation surfaces. 
Let $S$ be an oriented translation surface. Using the canonical charts for $S$ we can define one forms $dx$ and $dy$ on $S$. If $\sigma$ is a path in $S$ then we define the
holonomy vector corresponding to $\sigma$ to be the vector 
$$\biggl(\int_\sigma dx, \int_\sigma dy\biggr).$$
The holonomy of a path depends only on its homology class. Let $\Sigma$ be the set of singular points.  The holonomy gives a $\QQ$-linear map from  $H_1(S,\Sigma; \QQ)$ to $\RR^2$. We will call this map $\hol$. We call the image of this map the {\em holonomy} of $S$. In the language of Abelian differentials $\hol$ is called the period map. We can also consider the restriction of this map to $H_1(S;\QQ)$. We call the image of this map the {\em absolute holonomy} of $S$.

\begin{definition} 
\label{def:hd}
The set of {\em homological directions} is the set of directions of non-zero vectors in the holonomy.
\end{definition}

We will give a proposition which is useful in calculating the holonomy of a surface.

\begin{proposition}
\label{prop:hol criterion}
 If a translation surface $S$ can be obtained by gluing together polygons in the plane by gluing maps which are the restrictions of translations and all the vertices of these polygons lie in a subgroup $\Lambda\subset\RR^2$ then the holonomy of $S$ is contained in $\Lambda$.
 \end{proposition}
 
 \begin{proof} The surface $S$ decomposes into polygons. Any curve $\gamma$ in $S$ that begins and ends at singular points is homotopic rel endpoints to a sequence of edges of polygons $\gamma_j$. The holonomy of $\gamma$ is the sum of the holonomy of segments $\gamma_j$ and the holonomy of each segment lies in $\Lambda$.
 \end{proof}

The following definition identifies a type of cohomology class closely connected with the geometry of a translation surface. 

\begin{definition} A one dimensional cohomology class $\eta$ in $H^1(S,\Sigma; \RR)$ is {\em holonomic} if there exists a $Q$-linear function $\tau:\RR^2 \to \RR$ such that $\eta=\tau \circ \hol $.
\end{definition}

We have noted that the forms $dx$ and $dy$ are well defined one-forms on the translation surface $S$. The cohomology classes that they represent are examples of holonomic cohomology classes.
  If $\pi_1:\RR^2\to\RR$ denotes projection onto the first coordinate then the one form $dx$ represents the class corresponding to $\pi_1$. We call a linear combination of $dx$ and $dy$ a  {\em tautological one-form}.
  
The following is an alternative description of a holonomic cohomology class.  A one dimensional cohomology class  $\eta$ on a surface $S$ is holonomic if it corresponds to a $\QQ$-linear map from  $H_1(S; \QQ)$ to $\RR$ that vanishes on $ker(h)\cap H_1(S;\QQ)$. For a generic subspace of $H^1(S; \RR)$ the intersection with 
 $H^1(S; \QQ)$ will be $\{0\}$. It follows that for a generic translation surface every one dimensional cohomology class is holonomic.
 
We denote by $H^1(\RR^2;\RR)$ the group cohomology of $\RR^2$ considered as a discrete group. We have $H^1(\RR^2;\RR)=\Hom_\QQ(\RR^2,\RR)$.

The map:
$$H_1(S;\RR) \xrightarrow{\hol} H_1(\RR^2;\RR)$$
gives us
$$H^1(\RR^2;\RR) \xrightarrow{\hol^*} H^1(S;\RR).$$

\begin{proposition} The holonomic cohomology classes are those elements of $H^1(S;\RR)$ which are in the image of the map $\hol^*$.
\end{proposition}

The $J$ invariant can be approached very concretely as in the beginning of this section or rather abstractly and topologically. Here we give a topological interpretation of $J$.
The space $B\RR^2$ is a $K(\pi,1)$ for the group $\pi=\RR^2$ where $\RR^2$ is viewed as a discrete group. We can build the space $B\RR^2$ in two steps. First we build the simplicial space $E\RR^2$. The zero simplices of $E\RR^2$ are the points in $\RR^2$. The one simplices of $E\RR^2$ are pairs of points in $\RR^2$ and in general the n-simplices are $n+1$-tuples of points in $\RR^2$. The group $\RR^2$ acts simplicially on $E\RR^2$ and the quotient space is $B\RR^2$.

The cohomology of the space $B\RR^2$ is exactly the group cohomology of $\RR^2$. In particular $H_1(B\RR^2;\QQ)=\Hom(\RR^2;\QQ)$. We also note that $H^2(B\RR^2;\QQ)=\RR^2\wedge_\QQ\RR^2$. This gives us an alternative interpretation of $J(S)$.

\begin{proposition} If we identify $H^2(B\RR^2;\QQ)$ with $\RR^2\wedge_\QQ\RR^2$ then $J(S)=\hol_*([S])$.
\end{proposition}

\begin{proof} A triangulation of $S$ with geodesic edges gives an explicit map from $S$ to $B\RR^2$. A triangulation of $S$ lets us write $S$ as a collection of disjoint triangles $\Delta_j\subset\RR^2$ with identifications of edges. There is a map $\cup_j\Delta_j$ to $E\RR^2$ which takes each triangle to the corresponding two simplex in $E\RR^2$. If we compose with the quotient map from $E\RR^2$ to $B\RR^2$ then the resulting map has the property that points on parallel edges map to the same point. Thus the map from $\cup_j\Delta_j$ to $B\RR^2$ respects the identifications used to construct $S$ so it induces a map from $S$ to $B\RR^2$. The map induced on fundamental groups is the holonomy map and the image of the fundamental class under this map is $J$.
\end{proof}

The pairing on $H^1(S)$ which takes the one dimensional cohomology classes  $\alpha$ and $\beta$ to $\alpha\cup\beta$ evaluated on the fundamental class of $S$ is a skew-symmetric bilinear from which we denote by $\langle\alpha,\beta\rangle$. If $\alpha$ and $\beta$ are represented by differential forms we can write:

\begin{equation}
\label{eqn:skewform}
\langle\alpha,\beta\rangle=\int_S \alpha\wedge\beta.
\end{equation}

The $J$ invariant gives us a formula for evaluating this form on holonomic cohomology classes. If $\theta$ is a $\QQ$-linear map from $\RR^2$ to $\RR$ then define $\theta^*$ to be the cohomology class in $H^1(S;\RR)$ which corresponds to the homomorphism $\theta\circ\hol: H_1(S)\to\RR$.

 \begin{theorem} \label{thm: intJ} Let $\theta$ and $\tau$ be $\QQ$ linear maps from $\RR^2\to\RR$. Say that we can write $J(S)$ as $\sum_j v_j\wedge w_j$ then
 \begin{equation}
 \label{eqn: intJ}
\langle \theta^*, \tau^*\rangle=\frac{1}{2}\sum_j \theta(v_j) \tau(w_j)-\theta(w_j)\tau(v_j)
\end{equation}
 \end{theorem}
 
 \begin{proof} The expression on the right hand side is skew symmetric and $\QQ$ bilinear so it is independent of the particular way in which we write $J$.  Let us choose a simplicial triangulation of $S$ and with triangles $\Delta_j$. Let $v_j$ and $w_j$ be coordinates successive edges of $\Delta_j$. Then $J=\sum_j v_j\wedge w_j$. It suffices to prove the theorem for this particular expression for  $J$.
 
  It will be convenient to represent one dimensional cohomology classes as one forms which have certain discontinuities along boundaries of triangles.  According to (\cite{swan}) we can represent the real cohomology of a simplicial complex  in terms of piecewise polynomial forms on the simplicies of this complex which have certain compatibility relations on lower dimensional faces.  We can use these forms to calculate the cup product. One advantage of this representation is that the cup product operation is commutative on the level of cochains, not just on the level of cohomology, and this will be useful for us.

We can choose an orientation preserving affine parametrization of each triangle. We represent $\theta$ by a one form which is constant on each simplex and assigns the appropriate values to the edges of the simplex.

The restriction of $\theta$ to $\Delta_j$ is $\theta(v_j)dx+\theta(w_j)dy$ and the restriction of $\tau$ to $\Delta_j$ is $\tau(v_j)dx+\tau(w_j)dy$.

\begin{align*}
\int_{\Delta_j} \theta\wedge \tau&=\int_{\Delta} \bigl(\theta(v_j)dx+\theta(w_j)dy\bigr)\wedge \bigl(\tau(v_j)dx+\tau(w_j)dy\bigr)\\
&=\int_{\Delta} \theta(v_j)dx\wedge \tau(w_j)dy+\theta(w_j)dy\wedge\tau(v_j)dx\\
&=\int_{\Delta}\bigl( \theta(v_j) \tau(w_j)-\theta(w_j)\tau(v_j)\bigr)dx\wedge dy\\
&=\bigl( \theta(v_j) \tau(w_j)-\theta(w_j)\tau(v_j)\bigr)\int_{\Delta} dx \wedge dy \\
&=\frac{1}{2} \bigl( \theta(v_j) \tau(w_j)-\theta(w_j)\tau(v_j)\bigr)
\end{align*}

We can calculate the wedge product of the forms by summing over the triangles $\Delta_j$.

\begin{equation*}
\begin{split}
\int_S \theta\wedge \tau&=\sum_j \int_{\Delta_j} \theta\wedge \tau\\
&=\frac{1}{2}\sum_j \theta(v_j) \tau(w_j)-\theta(w_j)\tau(v_j)
\end{split}
\end{equation*}

\end{proof}

We can associate to a line $\ell\subset\RR^2$ the collection of holonomic cohomology classes corresponding to homomorphisms from $\RR^2$ to $\QQ$ with $\ell$ in their kernel.

\begin{proposition}[Levitt (see \cite{arnoux1})]
\label{prop:lagrangian}
The line $\ell$  represents an algebraically periodic direction exactly when this collection of cohomology classes is an isotropic subspace of the holonomic cohomology which is to say that the restriction of the cup product form to this subspace is zero. 
\end{proposition}

\begin{proof} Let $\pi_{\ell}: \RR^2 \to \RR$ be the quotient map with kernel $\ell$. (Here, we identify $\RR$ with $\RR^2/{\ell}$). The map $\pi_{\ell}$ induces a map $\hat{\pi}_{\ell}: \RR^2 \wedge_{\QQ} \RR^2 \to \RR \wedge_{\QQ} \RR$. Now, $\hat{\pi}_{\ell}(J)=\operatorname{SAF}_v(J)$ where $v$ is the direction of the line $\ell$. But the direction $v$ is algebraically periodic precisely when $\operatorname{SAF}_v(J)=0$.  Given homomorphisms $\sigma$ and $\tau$ from $\RR$ to $\QQ$, let $\sigma^*=\sigma \circ \pi_{\ell}$ and $\tau^*= \tau \circ \pi_{\ell}$ be the corresponding cohomology classes. If $\hat \pi_{\ell}(J)=0$, then $(\sigma \wedge \tau )(\hat \pi_{\ell}(J))=0$. Thus $\sigma^* \wedge \tau^*(J)=0$.
  The direction $v$ represents an algebraically periodic direction on $S$ precisely when $\operatorname{SAF}_v(J(S))=0$. 
\end{proof}

\end{section}

\begin{section}{the essential holonomy}

For lattice surfaces the holonomy field plays an important role. In the case of lattice surfaces this is an algebraic field related to the holonomy of $S$. In this section we will define a subspace of the absolute holonomy which we call the essential holonomy. We show that the essential holonomy is a subspace of $\QQ$-dimension at least 2 and no greater than the genus of $S$. While neither the holonomy nor the absolute holonomy are  scissors congruence invariants the essential holonomy is as we will show in Section~\ref{sec:brown}. In the next section we will describe the essential holonomy from a different point of view and show that it is the minimal $\QQ$ vector space of $\RR^2$ containing enough vectors to describe $J$.

We have denoted the pairing which takes the cohomology classes $\theta$ and $\tau$ to $\int_S\theta\wedge\tau$ by $\langle\theta,\tau\rangle$. This skew-symmetric bilinear form is non-degenerate on $H^1(S;\RR)$ but may become degenerate when restricted to the subspace of holonomic forms. Let $Z$ be the kernel of the restriction of this form.
Define $EC$ to be the quotient of the holonomic cohomology by $Z$.
 It follows that the pairing $\langle\cdot,\cdot\rangle$ is well defined on $EC$ and is non-degenerate. 

\begin{proposition} For a translation surface the dimension of $EC$ is an even number which is bounded below by two and above by twice the genus of $S$.
\end{proposition}

\begin{proof} Since $EC$ supports a non-degenerate skew symmetric bilinear form its dimension is even. Since $EC$ is a quotient of the group $H^1(S ; \RR)$ its dimension is bounded above by the dimension of $H^1(S ; \RR)$ which is $2g$. The forms $dx$ and $dy$ represent elements of $EC$. We have the fact that 
$$\langle dx, dy\rangle=\int_S dx\wedge dy=area(S)>0.$$
This implies that $dx$ and $dy$ are non-zero and represent linearly independent elements of $EC$. Thus the dimension is at least two.
\end{proof}

An element of $Z$ is a one dimensional cohomology class and it  can be evaluated on a one dimensional homology class.Thus each element of $Z$ defines a linear function from $H_1(S;\QQ)$ to $\RR$. We define the essential homology, $EHom$, to be the intersection of the kernels of these functions. We observe that two cohomology classes are equivalent if and only if they agree as linear functions on $EHom$ so we can identify $EC$ as $Hom_\QQ(EHom,\RR)$. Thus we have:

\begin{proposition} The dimension of $EHom$ as a rational vector space is even and bounded below by 2 and above by 2g.
\end{proposition}

\begin{proposition} $EHom$ is a subspace of the absolute homology of $S$.
\label{prop:ehominah}
\end{proposition}

\begin{proof} Let $\alpha\in H_1(S,\Sigma)$ be a homology class which is not in the absolute homology. Because $\alpha$ is not in the absolute homology we can find a cohomology class $\theta\in H^1(S,\Sigma)$ which vanishes on the absolute homology but which does not vanish when paired with $\alpha$. Now $\theta$ has the property that $\langle \theta, \tau\rangle=0$ for any $\tau$ since the value of $\langle \theta, \tau\rangle$ depends only on the values that $\theta$ takes on absolute homology classes.
In particular $\langle \theta, \tau\rangle=0$ for any holonomic cohomology class $\tau$ which is to say that $\theta\in Z$. On the other hand $\alpha$ is not in the kernel of $\theta$ so $\alpha$ is not in the absolute holonomy.
\end{proof}

\begin{definition} The {\em essential holonomy} is the rational subspace of $\RR^2$ which is  the image of $EHom$ under the holonomy map. 
\end{definition}

We denote the essential holonomy by $EH$. The holonomy map induces a bijection between $EHom$ and $EH$.

\begin{definition}
The  {\em absolute holonomy} of $S$ is the image of the of the absolute homology under the holonomy map. We denote this by $AH$.
\end{definition}

\begin{proposition} 
\label{prop:ehinah}
The essential holonomy $EH$ is contained in the asbolute holonomy $AH$.
\end{proposition}

\begin{proof} This follows from Proposition \ref{prop:ehominah} by applying the holonomy map.
\end{proof}

\end{section}

\begin{section}{A rational vector space determined by $J$}
\label{sec:brown}
In this section we show that $\Iso(J)$, the subgroup of elements of $\slr$ that preserve $J$, also preserves the essential holonomy. As we shall see, the essential holonomy may be characterized as the minimal $\QQ$-vector space $V_J$ such that $J$ can be written with entries in $V_J$, that is,
$$J= \sum_i v_i \wedge w_i$$ where $v_i, w_i \in V_J$. See Proposition~\ref{proposition:ehmin}.

The following argument was pointed out to us by Ken Brown.

\begin{proposition} \label{prop:brown}
Given an element $\beta \in \RR^2\wedge_\QQ\RR^2$, there exists a
 unique minimal  $\QQ$-vector space $V_{\beta} \subset \RR^2$ and an
element $\alpha \in H_2(V_{\beta};\QQ)$ so that if $\iota_*:
V_{\beta}\wedge V_\beta \to\RR^2\wedge_\QQ\RR^2$ is the map induced by the
inclusion $\iota: V_{\beta} \to \RR^2$, then $\iota_*(\alpha)= \beta$. Furthermore, for any $f \in GL(2,\RR)$ such that $f(\beta)=\beta$, $f(V_{\beta})=V_{\beta}$. \end{proposition}

\begin{proof} Consider the collection $\mathcal V$ of $\QQ$-vector subspaces of $\RR^2$ for which there is an $\alpha \in V\wedge_\QQ V$ with $\iota_{*}(\alpha)=\beta$.

First, note that there exists a finite dimensional $\QQ$-vector
space $V$ in $\mathcal V$.  We can write $\beta \in
\RR^2\wedge_\QQ\RR^2$ as $\beta=\sum_{i} v_i \wedge_{\QQ} w_i$. Let $V$
be the $\QQ$-vector space generated by the vectors $v_i$ and
$w_i$. Then $V$ has finite $\QQ$-dimension, and  there exists an
$\alpha \in V\wedge_\QQ V$ such that $\iota_*(\alpha)=\beta$.

$\mathcal V$ is closed under the operation of
taking intersections since for $V , W$ in $\mathcal V$,
$\bigwedge^2(V \cap W)=\bigwedge^2(V) \cap \bigwedge^2(W)$.
This implies that $\mathcal V$ contains
an element of minimal dimension. To prove uniqueness, let $V$ and $W$ be two
elements of $\mathcal V$ of minimal dimension then $V\cap W$ is
also in $\mathcal V$ and so it has  the same dimension as $V$ and
$W$. It follows that $V=W$.

Suppose $f \in GL(2,\RR)$ is such that $f(\beta)=\beta$. Then the
canonical vector space $V_{f(\beta)}=V_{\beta}$. Thus
$f(V_{\beta})=V_{\beta}$.

\end{proof}

\begin{corollary} \label{cor:isoj} Let $J \in \RR^2\wedge_\QQ\RR^2$. Then $\Iso(J)$ preserves the finite dimensional $\QQ$-vector space $V_J$ of $\RR^2$. \end{corollary}

\begin{proof} Note that $\Iso(J) \subset SL(2, \RR) \subset Aut(\RR^2)$ and, by definition, for any $f \in \Iso(J)$, $f(J)=J$. Then the corollary follows directly from the above proposition with $\beta = J$ and $V_J$ the associated canonical vector space. \end{proof}

Our goal is to show that $V_J$ is just
$EH(S)$.

 Let $\theta: AH \to
\QQ$ be a $\QQ$-linear map. Then $\theta$ corresponds to a
holonomic $1$-form on $S$.

\begin{proposition} If $\theta \neq 0$ is a holonomic $1$-form corresponding to a map from $AH \to \QQ$ and $\langle\theta,\tau\rangle = 0$ for all holonomic $1$-forms $\tau$, then $J(S)$ can be written as $\sum_i v_i \wedge w_i$  where $v_i, w_i \in ker(\theta)$. 
\end{proposition}

\begin{proof}
Let $\lbrace e^0, e^1, \cdots, e^n \rbrace$ be a basis of $AH^*$
where $e^0=\theta$. Let $\lbrace e_0, \cdots , e_n \rbrace$ be the
dual basis for $AH$. Then $$\lbrace e_i \wedge e_j \rbrace_{0 \leq
i < j \leq n}$$ is a basis for $\bigwedge^2_{\QQ} AH$ and
$$
\lbrace e^i \wedge e^j \rbrace_{0 \leq i < j \leq n}
$$ is the
dual basis for $\bigwedge_{\QQ}^2 AH^*$. As we noted above, since
$J=\sum_i v_i \wedge w_i$ for $v_i, w_i \in AH$, we have that
$$
J=\sum_{0 \leq i < j \leq n} r_{ij} e_i \wedge e_j.
$$ 
 Then
\begin{equation} \label{eq:4} 
\begin{split}
 \langle e^i,e^j\rangle &= \int_S e^i \wedge e^j \\
&=\sum_{0 \leq k < m \leq n} \langle e^i \wedge e^j, r_{km} e_k \wedge e_m\rangle \\
&=r_{ij}. \\ 
\end{split}
\end{equation}

By hypothesis, $\langle e^0, e^j \rangle = 0$ for all $j$ since $e^0 = \theta$. So $r_{0j}=0$ for all $j$. It follows that we can write
$J= \sum_{1 \leq i < j \leq n} r_{ij} e_i \wedge e_j$. But $e_i
\in ker(\theta)$ for $i \neq 0$.
\end{proof}

Recall that $Z$ is the subset of holonomic cohomology consisting of
those forms $\theta$ such that $\langle\theta, \tau\rangle=0$ for every
holonomic one-form $\tau$ and that the  essential holonomy of S, $EH(S)$, is the subspace of $im(h)$ on which $Z$ vanishes.

\begin{proposition} We can write $J= \sum_i v_i \wedge w_i$ where $v_i, w_i \in EH$. \end{proposition}

\begin{proof}
Since $\Ann(Z)=EH$, by repeatedly applying the previous proposition
for each $\theta \in Z$, we find that $J=\sum v_i \wedge w_i$
where $v_i, w_i \in EH$.
\end{proof}

\begin{proposition} \label{proposition:ehmin} $EH$ is the minimal $\QQ$-vector space $V$ such that $J=\sum_i v_i \wedge w_i$ where $v_i, w_i \in V$. \end{proposition}
\begin{proof} Suppose for a contradiction that there exists a $\QQ$-vector space $V$ so that $V$ is strictly contained in $EH$. Then there exists a non-zero $\QQ$-linear map $\theta : EH \to \QQ$ so that $V \subset ker(\theta)$. Then $\theta$ represents a holonomic one-form which we  call $\theta^*$. Then for any other form $\tau$,  $\langle \theta^* , \tau \rangle = 0$. Since the form is non-degenerate this implies that $\theta^*=0$ contradicting our assumption.
\end{proof}

\begin{corollary} \label{corollary:v_j=eh} If $J(S) \in \RR^2\wedge_\QQ\RR^2$ then $V_J = EH(S)$. \end{corollary}
\begin{proof} This follows from Proposition~\ref{proposition:ehmin} and the definition of $V_J$. \end{proof}

We take advantage of the fact that the essential holonomy of $S$ depends only on $J(S)$ and not on the surface $S$ to make the following definition.

\begin{definition} Given $J\in \RR^2\wedge_\QQ\RR^2$ we define $EH(J)$ to be $V_J$.
\end{definition}

Corollary~\ref{corollary:v_j=eh} shows that this new definition is consistent with our previous definition.

Since the $EH$ depends only on $J$ by Corollary~\ref{corollary:v_j=eh} and since $J$ is a scissors congruence invariant, we have that 

\begin{theorem} The essential holonomy is a scissors congruence invariant.
\end{theorem}
On the other hand, in Theorem \ref{eh ne ah} we show that the holonomy and even the genus of a translation structure are not scissors congruence invariants.

We note that the construction of $V_J$ is of a general nature and does not depend on the fact that $J$ comes from a translation surface $S$. Thus we can extend the definitions of certain objects associated with translation surfaces to the more abstract context of arbitrary elements $J\in \RR^2\wedge_\QQ\RR^2$. These include the area, the essential holonomy and the set of algebraically periodic directions. This observation will be useful in Section~ \ref{sec:construction} when we want to construct translation surfaces with prescribed periodic direction fields. In this case we will first construct an abstract  $J$ with the desired properties and then we will construct a surface that realizes this $J$.

\begin{corollary} \label{corollary:JEH} $\Iso(J)$ preserves $EH(S)$. \end{corollary}
\begin{proof} This follows directly from Corollary~\ref{corollary:v_j=eh} and Proposition~\ref{cor:isoj}.
\end{proof}

Next, we will give an alternate interpretation of $\Iso(J)$.
Recall that there is an action of $GL(2,\RR)$ on $1$-forms so that if $g \in GL(2,\RR)$ and $\eta$ is a $1$-form, then $g\cdot \eta=\eta \circ g$. Note that this action is natural with respect to the action of $GL(2,\RR)$ on surfaces: if $S$ is a surface, then the translation structure defined on $gS$ is given by precomposition of $g$ with coordinate charts of $S$.

We say that $\alpha\in\slr$ preserves the product pairing if
\[\int_S \alpha(\sigma)\wedge\alpha(\tau)=\int_S \sigma\wedge\tau
\]
for all pairs $\sigma$ and $\tau$.

\begin{theorem} 
\label{thm:skewform}The skew symmetric form $\langle\cdot,\cdot\rangle$ on the essential holonomy determines $J$.
\end{theorem}
\begin{proof}  Let $e_1\ldots e_n$ be a $\QQ$-basis for $EH$. We can write $J$ as 
$$
J=\sum_{0 \leq i < j \leq n} r_{ij} e_i \wedge e_j.
$$ 
Let $e^0\ldots e^n$ be  a dual basis. Equation~\ref{eq:4} shows that $r_{ij}=\langle e^i,e^j\rangle$. 
\end{proof}

\begin{corollary} \label{thm:HSG} $\Iso(J(S))$ is exactly the collection of elements of $\slr$ that preserve the product pairing when acting on holonomic 1-forms.
\end{corollary}

\begin{proof} Theorem~\ref{thm: intJ} shows that $J$ determines the form $\langle\cdot,\cdot\rangle$ on ${\RR^2}^*$. The previous theorem shows that $\langle\cdot,\cdot\rangle$ determines $J$. Thus any linear map of $\RR^2$ that preserves $\langle\cdot,\cdot\rangle$ preserves $J$.
\end{proof}

\begin{proposition} If $S$ and $S'$ are scissors congruent then $EH(S)=EH(S')$ and the bilinear forms induced on $EH^*$ for both surfaces are the same.
\end{proposition}

\begin{proof} Scissors congruence preserves the $J$ invariant. The group $EH$ depends only on the $J$ invariant. The bilinear form on $EH$ is determined by $J$.
\end{proof}

\end{section}

\begin{section}{The structure of $EH$} \label{sec:struceh}

In this section, we examine $EH(S)$ as a rational subspace of $\RR^2$. Since the results are of a general algebraic nature we state them for a general finite dimensional rational subspace of $\RR^2$ which we denote by $\Lambda$ rather than stating them for $EH(S)$. In the first half of this section we do not use the skew-symmetric form on $EH(S)$ which arises from $J$. 
We continue the analysis by assuming that  $\Lambda=V_J$ for some element $J\in\RR^2\wedge_\QQ\RR^2$ as in  section~\ref{prop:brown}.

Let $\Lambda$ be a finite dimensional rational subspace of $\RR^2$ which has the property that it contains a pair of vectors that generate $\RR^2$ as a real vector space. In case $\Lambda=EH(s)$ we have $\dim_\QQ(\Lambda)\le 2g$ where $g$ is the genus of $S$. 

\begin{definition} The {\it multiplicative field} of $\Lambda$ is defined to be $k=\lbrace \lambda \in \RR : \lambda \Lambda \subset \Lambda \rbrace$. In the case that $\Lambda=EH(S)$ we will call $k$ the multiplicative field of $S$. 
\end{definition}

\begin{theorem}\label{thm:kalg}
The field $k$ is algebraic and the degree of this field over $\QQ$ is bounded above by $\dim_\QQ(\Lambda)/2$. 
\end{theorem}

\begin{proof} By assumption $\Lambda$ contains vectors $v_1$ and $v_2$ which are linearly independent (over $\RR$). Let $\ell_j$ be the one dimensional real subspaces spanned by the$v_j$. Let $r_j$ be the $\QQ$-dimension of $\ell_j\cap\Lambda$. Since $\Lambda\cap\ell_j$ is a $k$-vector space of dimension at least one we have $\dim_\QQ(k)\le r_j$. In particular $k$ is a finite extension of $\QQ$. 

The subspace of $\Lambda$ generated by $\ell_1\cap\Lambda$ and $\ell_2\cap\Lambda$ has dimension $r_1+r_2$. Thus $r_1+r_2\le\dim_\QQ(\Lambda)$. Combining this with the previous inequality we have $2\dim_\QQ(k)\le\dim_\QQ(\Lambda)$.
\end{proof}

\begin{corollary} The multiplicative field of $S$ is algebraic and its degree is bounded by the genus of $S$.
\end{corollary}

\begin{proof} This follows because the dimension of $EH(S)$ is bounded by twice the genus of $S$.
\end{proof}

\begin{definition} A line $\ell\subset\RR^2$  is {\it maximal} if $\dim_\QQ(\Lambda\cap\ell)\ge\dim_\QQ(\Lambda)/2$.
\end{definition}

Next we will consider those $\QQ$ vector spaces $\Lambda$ with many maximal lines. 

\begin{proposition} \label{prop:k}
If there are three maximal lines in $\Lambda$  there are infinitely many. If the lines of slope slopes $0, 1$ and $\infty$ are maximal then $\Lambda$  then the maximal lines are those with slopes in $k\cup\infty$ where $k$ is the multiplicative field of $\Lambda$.
Furthermore $\Lambda=\{(x,y): x,y\in\Gamma\}$ where $\Gamma\subset\RR$ is a finite dimensional $k$ vector space.
\end{proposition}

We remark that if there are three maximal lines we may choose coordinates so that they are the lines with slopes $0, 1$ and $\infty$.

\begin{proof} To begin assume that  there are two distinct maximal lines $\ell_1$ and $\ell_2$. Let $r_j$ be the $\QQ$-dimension of $\ell_j\cap\Lambda$ so $r_j\ge r$. The subspace of $\Lambda$ generated by $\ell_1\cap\Lambda$ and $\ell_2\cap\Lambda$ has dimension $r_1+r_2$. Thus $r_1+r_2\le \dim_\QQ(\Lambda)$. Since $r_j\ge\dim_\QQ(\Lambda)/2$ by assumption we conclude that $r_1=r_2=\dim_\QQ(\Lambda)/2$ and furthermore
$\Lambda$  is generated by  $\Lambda\cap\ell_1$ and $\Lambda\cap\ell_2$.

Now choose coordinates so that $\ell_1$ and $\ell_2$ are the standard coordinate axes. Thus $\Lambda$ is the set of pairs $(x,y)$ such that $x\in\ell_1\cap \Lambda=\Lambda_x$ and $y\in\ell_2\cap \Lambda=\Lambda_y$.

Consider a third line $\ell'$ with slope $\lambda$ which we do not assume to be maximal. A point $(x,y)$ is in $\ell'$ if and only if $y=\lambda x$ where
$x\in\Lambda_x$ and $y\in\Lambda_y$. Thus the dimension of $\Lambda\cap\ell'$ is the dimension of $\lambda\Lambda_x\cap\Lambda_y$. In particular by looking at dimensions we see that $\ell'$ is maximal if and only if $\lambda\Lambda_x=\Lambda_y$. 
Now assume that the line with slope 1 is maximal. This means that $\Lambda_x=\Lambda_y$.

Let $k=\{\lambda\in\RR: \lambda \Lambda_x \subset \Lambda_x\}$. The maximal lines are just those with slopes in $k\cup\infty$. Now let $\Gamma=\Lambda_x$.
\end{proof}

We defined the multiplicative field to be the set $k=\{\lambda\in\RR: \lambda \Lambda \subset \Lambda\}$
As a consequence of the above proof we see that if the lines of slopes $0,1$ and $\infty$ in $\Lambda$ are maximal then the multiplicative field can also be described as $k=\{\lambda\in\RR: \lambda \Gamma \subset \Gamma\}$. In particular $\Gamma$ is a finite dimensional $k$ vector space.

We see that $\Lambda$ can be written as $\Lambda=\alpha_1k^2 \oplus \cdots \oplus \alpha_n k^2$ where $n$ is half the $k$ dimension of $\Lambda$ and $\alpha_1\ldots\alpha_n$ are real numbers which give a $k$ basis for $\Gamma$. In particular we see that the dimension of $\Lambda$ as a $k$ vector space is even.
If the $k$ dimension of $\Lambda$ is two then  by rescaling we can write $\Lambda=k^2$.

\begin{proposition} Every non-zero vector in $\Lambda$ determines a maximal direction if and only if the $k$ dimension of $\Lambda$ is two.
\label{prop:rank1}
\end{proposition}
\begin{proof} If the $k$ dimension of $\Lambda$ is two then a maximal line is a line $\ell$ so that $\dim_\QQ(\Lambda\cap\ell)\ge\dim_\QQ(k)$. If $\ell$ contains a non-zero vector of $\Lambda$ then it contains a one dimensional $k$ vector space so this condition holds.

Now we prove the other direction.
The hypothesis that  every line $\ell$ for which $\Lambda\cap\ell\ne0$ is maximal implies that
there are infinitely many maximal lines.  Thus we may assume that the lines of slopes $0,1$ and $\infty$ in $\Lambda$ are maximal. Consider a line $\ell'$ of slope $\lambda$. To say that $\ell'\cap\Lambda\ne0$ is to say that
 $\lambda\Lambda_x\cap\Lambda_x\ne\{0\}$.   To say that $\ell'$ is maximal is to say that $\lambda\Lambda_x=\Lambda_x$. The hypothesis that the first condition implies the second means that  multiplication by non-zero elements of $k$ acts transitively on the non-zero elements of $\Lambda_x$. This implies that $\Lambda_x$ must be a one-dimensional vector space over $k$. Thus according to the previous theorem we may assume that $\Lambda=\alpha k^2$. By a further coordinate change we may assume that $\alpha=1$.
\end{proof}

\begin{lemma} Assume that $\Lambda$ has three maximal directions and and that their slopes are $0, 1$ and $\infty$. Then the subgroup of $GL(2,\RR)$ that preserves $\Lambda$ is  $GL(2,k)$ where $k$ is the multiplicative field of $\Lambda$.
\end{lemma}
\begin{proof}  We have 
$$
\Lambda=\alpha_1k^2 \oplus \cdots \oplus \alpha_n k^2.
$$ 
Thus $GL(2,k)$ is contained in the affine automorphism group of $\Lambda$ since it is a sum of vector spaces invariant under $GL(2,k)$. 

Say we have a linear transformation $M$ that preserves $\Lambda$. Then it also preserves the collection of maximal directions which are just the lines with slopes in $k$. 
Thus $M$ takes the triple of directions $(0,1,\infty)$ to another triple of distinct directions $(s_1, s_2, s_3)$ where each $s_i \in k$.
The group $GL(2,k)$ acts transitively on  triples of distinct lines with slopes in $k$, thus there exists a matrix $M' \in GL(2,k)$ which takes $(0,1,\infty)$ to $(s_1, s_2, s_3)$. In particular the matrix $(M')^{-1} M$ takes the triple $(0,1,\infty)$ to itself. This implies that $(M')^{-1} M = \lambda I$ where $\lambda I$ preserves $\Lambda$ so $\lambda \in k$.  Thus $M=\lambda M'\in GL(2,k)$. This completes the proof.
 \end{proof}

 For the remainder of this section we will put some additional assumptions on $\Lambda$.
We assume that we have some $J\in \RR^2\wedge_\QQ\RR^2$ and that $\Lambda=V_J$. The element $J$ defines a skew-symmetric non-degenerate inner product $\langle,\rangle$ on $Hom(\Lambda,\QQ)$. As we showed in section 4 when $J=J(S)$ then $V_J=EH(S)$. The case when $J=J(S)$ is the case of primary interest though it will be convenient to deal with other elements of $\RR^2\wedge_\QQ\RR^2$ which are not a priori of this form.
 
Many properties that we have associated with a translation surface $S$ depend only on $J(S)$. We take advantage of this fact to define the corresponding properties abstractly for an arbitrary element  $J\in\RR^2\wedge_\QQ\RR^2$.

Let $\ell$ be a one dimensional real vector subspace of $\RR^2$. Let $\Lambda_\ell^*$ be the subspace of  $Hom(\Lambda,\QQ)$ consisting of elements that vanish on $\ell\cap\Lambda$.

We say that a line $\ell$ is {\em algebraically periodic} if $\langle \theta,\tau\rangle=0$ for any $\theta$, $\tau$ in $\Lambda_\ell^*$. If we choose $J$ to be $J(S)$ for some surface $S$ then this notion of an algebraically periodic direction coincides with the previous notion. 
 
\begin{definition} $J\in\RR^2\wedge_\QQ\RR^2$ is algebraically periodic if there exist at least three directions $v$ such that $SAF_v(J)=0$. \end{definition}

\begin{definition} We say that $J$  is in standard form if the $SAF$-invariant vanishes along the lines with slopes $0, 1$ and $\infty$. \end{definition} 

 If $S$ is in standard form, then a theorem from (\cite{calta}) implies that there are infinitely many algebraically periodic directions which form the projective line $\PP(K)$ where $K$ is a field. The theorem is proved by analyzing the collection of parabolic matrices that preserve $J$. Thus we have the following.

\begin{proposition} If $J$ is algebraically periodic and in standard form then the slopes of periodic directions form a number field $K$.
\end{proposition}

\begin{definition} Given $J\in\RR^2\wedge_\QQ\RR^2$ we define the field $K$ from the previous proposition to be the periodic direction field.
\end{definition}

The key point here is that any $J$ with three algebraically periodic directions can be put in standard form by means of a coordinate change.

\begin{definition} We say that a surface $S$ is in standard form if the lines with slopes $0, 1$ and $\infty$ are algebraically periodic. \end{definition}

 \begin{proposition} The lines $\ell$ in $\Lambda$ which represent algebraically periodic directions  are maximal.
 \end{proposition}

\begin{proof}  Let $\ell$ be an algebraically periodic line in $\RR^2$. Let $\Lambda_\ell^*$ be the subspace of  $Hom(\Lambda,\QQ)$ consisting of elements that vanish on $\ell\cap\Lambda$. Since the line $\ell$ is algebraically periodic $\langle \theta,\tau\rangle=0$ for any $\theta$, $\tau$ in $\Lambda_\ell^*$. This means that $\Lambda_\ell^*$ is an isotropic subspace for the non-degenerate skew-symmetric form $\langle\cdot,\cdot\rangle$. The dimension of an isotropic subspace is at most $\dim\Lambda/2$. On the other hand $\dim_\QQ(\Lambda_\ell^*)+\dim_\QQ(\Lambda\cap\ell)=\dim_\QQ(\Lambda)$ so the dimension of $\Lambda\cap\ell$ is at least $\dim\Lambda/2$ and $\ell$ is maximal.
\end{proof}

\begin{proposition} \label{prop:apk} If $J$ is algebraically periodic and in standard form, the periodic direction field $K$ is contained in the multiplicative field $k$ in particular $K$ is a number field.
\end{proposition}

\begin{proof}
The lines $\ell$ in $\Lambda$ whose slopes are algebraically periodic directions on $S$ are maximal. But Proposition~\ref{prop:k} shows that the maximal lines are those with slopes in $k \cup \lbrace \infty \rbrace$.
\end{proof}

The following definition is standard. 

\begin{definition} The {\it holonomy field} of $S$ is the smallest field $L$ such that the absolute holonomy is contained in a two-dimensional
$L$-vector space. \end{definition}

\begin{proposition} \label{prop: k,K}
The multiplicative field of $S$ is contained in the holonomy field.
 \end{proposition}

\begin{proof} Say $\lambda$ is a non-zero element of the multiplicative field $k$. Let $v$ be an element
of $EH$. The definition of $k$ gives us that $w=\lambda v$ is also
in $EH$. We use the fact that $AH$ generates $\RR^2$ as an $\RR$
vector space to find an element $u$ which is not an $\RR$ multiple
of $w$. Since $AH$ generates a two dimensional $L$ vector space
there must be a $L$ linear relation between
 $u$, $v$ and $w$. Thus $ru+sw+tv=0$ where the coefficients are in $L$ and are
 not all zero. Since $u$ and $v$ are $\RR$ linearly independent we have $r=0$.
 So $w=-(t/s)v$. Thus $\lambda=-(t/s)$ is an element of $L$.
\end{proof}

\begin{theorem} \label{theorem:edap} Every essential homological direction is algebraically periodic if and only if $EH$ has dimension two over the periodic direction field $K$.
\end{theorem}

\begin{proof} This follows from Proposition~\ref{prop:rank1} and the fact that every algebraically periodic direction is maximal.
\end{proof} 

\begin{definition} 
A surface $S$ is {\em completely algebraically periodic} if every homological direction is algebraically periodic. 
 \end{definition} 
 
 The following result gives a useful criterion for recognizing completely algebraically periodic surfaces.

\begin{corollary} A translation surface $S$ is completely algebraically periodic if it  can be built by gluing a collection of polygons each of which has vertices in the periodic direction field.
\end{corollary}

\begin{proof} This follows from Proposition~\ref{prop:hol criterion}.
\end{proof}

\begin{theorem} \label{theorem:ahcp} An algebraically periodic surface $S$ is completely algebraically periodic if and only if the holonomy is a vector space over the periodic direction field of dimension 2.
\end{theorem}

\begin{proof} We prove the reverse implication first. If the absolute holonomy is a vector space over the periodic direction field $K$ of dimension 2,
then since the absolute holonomy contains the essential holonomy and the essential holonomy is a vector space over $K$ of dimension at least 2, the absolute holonomy is equal to the essential holonomy. So the $K$-dimension of $EH$ is two and Theorem~\ref{theorem:edap} implies every essential homological direction is algebraically periodic and so every absolute homological direction is algebraically periodic. Thus $S$ is completely algebraically periodic.

To prove the forward implication, suppose that $S$ is completely algebraically periodic. Then every homological direction is algebraically periodic. By Theorem~\ref{theorem:edap} we may assume that $EH=k^2\subset\RR^2$.

From Proposition~\ref{prop:ehinah} we know that $EH\subset AH$. It suffices to show that $EH=AH$ for then Theorem~\ref{theorem:edap} can be used. Let $p\in AH$. Since the sets of directions for $EH$ and $AH$ are the same by Proposition~\ref{prop:rank1}, 
we know that $p$ lies in a line $\ell$ through the origin with slope in $K$ so $\ell$ is given by an equation with coefficients in $K$. Let $q$ be a point in $EH$ which does not lie on this line. Let $\ell'$ be the unique line through $p$ and $q$. The slope of the line $\ell'$ lies in $K$ and this line contains a point with coordinates in $K$ so this line is given by an equation with coefficients in $K$. Now $p=\ell\cap\ell'$ and since both lines are given by equations with coefficients in $K$ their intersection must have both coordinates in $K$. So $p\in EH$ as was to be shown.
\end{proof}

\begin{corollary}
\label{cor:capeq}
A surface is completely algebraically periodic if and only if its holonomy field is equal to its periodic direction field.
\end{corollary}

\begin{proof} Let $L$ be the holonomy field of $S$. If $S$ is completely algebraically periodic, then Theorem~\ref{theorem:ahcp} implies that $\dim_K(AH)=2$. But the definition of the holonomy field implies that $L \subset K$.  On the other hand, Proposition~\ref{prop: k,K} and Proposition $5.5$ imply that $K \subset L$. Thus, $L=K$. 
If $L=K$, then $\dim_K(AH)=2$ and Theorem~\ref{theorem:ahcp} implies that $S$ is completely algebraically periodic. \end{proof} 

In the section 5 we show that $EH(J)$ is $\Gamma\oplus\Gamma$ where $\Gamma$ is a $K$ vector space. It follows that $EH(J)$ is a $K$ vector space of even dimension.

\begin{definition} We define the rank of $J$ to be $\frac{1}{2}\dim_K(EH(J))$.
\end{definition}

\begin{proposition} If $S$ is completely algebraically periodic then $J(S)$ has rank one.
\end{proposition}

\begin{proof} This follows from Theorem~\ref{theorem:edap}.
\end{proof}

We end this section with a criterion for recognizing rank one $J$'s.

\begin{definition} We say that $J$ can be written with coefficients in a field $L$ if $J$ is in standard form and $J$ can be expressed as 
$$J= \sum_i v_i \wedge w_i$$ where the components of each $v_i$ and $w_i$ are elements of $L$. \end{definition}

The results of section 4 imply that this is equivalent to saying that $EH(J)=L^2$.

\begin{proposition}
\label{prop:characterization} If $J$ can be written with coefficients in a field $L$ and all lines with slopes in $L\cup\infty$ are algebraically periodic then $\Lambda_J=L^2$ and $J$ has rank one with periodic direction field equal to $L$.
\end{proposition}
\begin{proof} 
By Proposition~\ref{proposition:ehmin}, the assumptions imply that $\Lambda=\Lambda_J\subset L^2$ and that the periodic direction field contains $L$. Here $\Lambda_J$ is defined as in the statement of Proposition~\ref{prop:brown}. Let $k$ be the multiplicative field. Let $m$ be the degree of the multiplicative field and $d$ the degree of the direction field. We know that $K\subset k$ so $d\le m$. 
On the other hand $\Lambda$ is a vector space over $k$ of dimension at least 2 so $2m\le \dim_\QQ(\Lambda)$. $\Lambda\subset L^2$ so $\dim_\QQ(\Lambda)\le 2\deg(L)$. $L\subset K$ so $\deg(L)\le d$. Thus $2m\le2d$. We conclude that  $m=d$ so $L=K$ and $\Lambda=L^2=K^2$.
\end{proof}

\end{section}

\begin{section}{The Periodic Direction Field}
\label{sec:pdf}

In this section we describe a symmetric bilinear form on $\Gamma^*$ which is used in (\cite{kenyon-smillie}) and related to (\cite{kenyon}). 
When combined with Theorem~\ref{thm: intJ}, we get Theorem~\ref{thm:1.6-1} which gives criteria for algebraic periodicity of a surface $S$ in terms of explicit equations involving the coordinates of the polygons from which $S$ is constructed. 

As shown in section 5 we have $\Lambda=\Gamma\oplus\Gamma$ where $\Gamma\subset\RR$. We also have $\Lambda^*=\Gamma^*\oplus \Gamma^*$.

Here, we introduce notation that will be used throughout this section. The symbol 
$$\sigma dx\ (\rm{resp.\ } \sigma dy)$$ 
which lies in $\Lambda^*$ represents the homomorphism which takes a vector $\begin{bmatrix} x \\ y \end{bmatrix}$ in $\Lambda$ to the complex number $\sigma(x)\ (\rm{resp.\ } \sigma(y))$. 

\begin{theorem} 
\label{alg per crit} 
The surface $S$ is algebraically periodic in standard form if and only if for every $\sigma$ and $\tau$ in $\Gamma^*$ we have:
\begin{align}
\langle \sigma dx, \tau dx\rangle &=0\label{eq:dx}\\
\langle \sigma dy, \tau dy\rangle&=0\\
\langle \tau dx, \sigma dy\rangle&=\langle \sigma dx, \tau dy\rangle \label{eq:cross} 
\end{align}
\end{theorem}
\begin{proof}

Since the horizontal direction is algebraically periodic by Proposition~\ref{prop:lagrangian} we have $\langle \sigma dx, \tau dx\rangle=0$ for all $\sigma$ and $\tau$ in $\Gamma^*$.
Since the vertical direction is algebraically periodic we have $\langle \sigma dy, \tau dy\rangle=0$ for all $\sigma$ and $\tau$ in $\Gamma^*$.

Since the line with slope 1 is algebraically periodic we have 
$$\langle \sigma (dy-dx), \tau (dy-dx)\rangle=0$$ for all $\sigma$ and $\tau$ in $\Gamma^*$.
We expand the left hand side and use the observation about the horizontal and vertical directions.
\begin{equation*}
\begin{split}
\label{symmetry}
0&=\langle (\sigma dy-\sigma dx), (\tau dy-\tau dx)\rangle\\
&= \langle \sigma dy, \tau dy\rangle-\langle \sigma dy, \tau dx\rangle-\langle \sigma dx, \tau dy\rangle +\langle \sigma dx, \tau dx\rangle\\
&=-\langle \sigma dy, \tau dx\rangle-\langle \sigma dx, \tau dy\rangle\\
&=\langle \tau dx, \sigma dy\rangle-\langle \sigma dx, \tau dy\rangle.
\end{split}
\end{equation*}
Equation~\ref{eq:cross} follows. 
\end{proof}

\begin{definition} Let $\sigma$ and $\tau$ be elements of $\Gamma^*$. Then we define $$(\sigma,\tau)=\langle \sigma dx,\tau dy\rangle.$$ \end{definition}

This symmetric bilinear form is closely related to one that appears in the work of Kenyon (\cite{kenyon}).

\begin{corollary}
\label{cor:symform}
 If $S$ is algebraically periodic and in standard form then $(\cdot, \cdot)$ is a non-degenerate symmetric, bilinear form on $\Gamma^*$.
\end{corollary}
\begin{proof}
The symmetry follows from Equation~\ref{eq:cross}. To see that it is non-degenerate we use the fact that the form $\langle\cdot,\cdot\rangle$ is non-degenerate. If $\sigma$ is a non-zero element of $\Gamma^*$ then there is some element $\tau=\tau_1dx+\tau_2dy$ of $\Lambda^*$ such that $\langle\sigma dx,\tau\rangle\ne0$.  We have
\begin{align*}
\langle\sigma dx,\tau\rangle&=\langle\sigma dx,\tau_1dx+\tau_2dy\rangle\\
&=\langle\sigma dx,\tau_1dx\rangle+\langle\sigma dx,\tau_2dy\rangle\\
&=\langle\sigma dx,\tau_2dy\rangle\\
&=(\sigma,\tau_2)
\end{align*}
where the fact that $\langle\sigma dx,\tau_1dx\rangle=0$ follows from Equation~\ref{eq:dx}. 
Thus given $\sigma\ne0$ we have produced a $\tau_2$ so that $(\sigma,\tau_2)\ne0$.
\end{proof}

\begin{proposition}
\label{prop:innerprod}
 If $S$ is algebraically periodic and in standard form then the form $(\cdot,\cdot)$ determines $J$.
\end{proposition}

\begin{proof} Theorem \ref{thm:skewform} shows that $\langle\cdot,\cdot\rangle$ determines $J$. Here we show that $(\cdot,\cdot)$ determines $\langle\cdot,\cdot\rangle$. Let $\alpha$ and $\beta$ be  elements of $\Lambda^*$. We can write each of these in terms of its components as $\alpha=\alpha_x dx+\alpha_y dy$ and 
$\beta=\beta_x dx+\beta_y dy$. Then, since $S$ is algebraically periodic and in standard form, 

\begin{equation*}
\begin{split}
\langle\alpha,\beta\rangle&=\langle \alpha_x dx+\alpha_y dy, \beta_x dx+\beta_y dy\rangle\\
&=\langle \alpha_x dx , \beta_x dx  \rangle+
\langle\alpha_x dx  , \beta_y dy  \rangle+
\langle\alpha_y dy  ,  \beta_x dx  \rangle+
\langle\alpha_y dy  , \beta_y dy  \rangle\\
&=
\langle\alpha_x dx  , \beta_y dy  \rangle-
\langle  \beta_x dx ,  \alpha_y dy \rangle\\
&=(\alpha_x,\beta_y)-(\beta_x,\alpha_y).
\end{split}
\end{equation*}
\end{proof}

Recall that $k\subset\RR$ is the multiplicative field of $S$ and that $\Gamma$ is a finite dimensional $k$-vector space. We will find an explicit basis for $\Gamma^*$ and $\Lambda^*$.

Let $\phi^p:\Gamma\to k$ for $p=1 ,\ldots,n$ be a basis of $k$-linear maps.
Let $\phi_q: k\to\CC$ for $q=1\ldots m$ be the distinct field embeddings. Let $\phi^p_q:\Gamma\to \CC$ be defined as $\phi_q^p=\phi_q \circ \phi^p$. The collection of linear maps $\phi^p_q$ give a basis for $\Gamma^*$.
For $\gamma\in\Gamma$ and $\alpha\in k$ we have $\phi^p_q(\alpha\gamma)=\phi_q(\alpha)\phi^p_q(\gamma)$.
Let $dx$ and $dy$ represent the coordinate projections from $\RR^2$ to the $x$ and $y$ axes. 
The linear maps $\phi^p_q dx$ and $\phi^p_q dy$ give a basis for $\Lambda^*$.

\begin{proposition} The line with slope $\lambda$ is algebraically periodic if $\phi_s(\lambda)=\phi_q(\lambda)$ whenever $\langle\phi^r_s dx, \phi^p_q dy\rangle\ne 0$.
 \end{proposition}

\begin{proof} By Proposition~\ref{prop:apk} $K \subset k$. 
The form $dy-\lambda dx$ gives a projection from $\RR^2$ to $\RR$ for which the kernel is the line of slope $\lambda$. This line will correspond to an isotropic subspace if 

$$\langle\phi_q^p(dy-\lambda dx), \phi_s^r(dy-\lambda dx)\rangle=0$$

for all $p,q,r,s$. 

Expanding this expression and using Equation~\ref{eq:cross} we get:

\begin{align*}
0&=\langle \bigl(\phi^p_q dy-\phi_q(\lambda)\phi^p_q dx\bigr),\bigl(\phi^r_s dy-\phi_s(\lambda)\phi^r_s dx\bigr)\rangle\\
&=\langle -\phi_s(\lambda)\phi^p_q dy, \phi^r_s dx\rangle-\langle\phi_q(\lambda)\phi^p_q dx, \phi^r_s dy\rangle\\
&=\langle -\phi_s(\lambda)\phi^r_s dy, \phi^p_q dx\rangle-\langle\phi_q(\lambda)\phi^p_q dx, \phi^r_s dy\rangle\\
&=\langle\phi_s(\lambda)\phi^p_q dx, \phi^r_s dy\rangle-\langle\phi_q(\lambda)\phi^p_q dx, \phi^r_s dy\rangle\\
&=(\phi_s(\lambda) -\phi_q(\lambda))\langle\phi^r_s dx, \phi^p_q dy\rangle. 
 \end{align*}
\end{proof}

\begin{corollary} 
\label{cor:k=K}
The multiplicative field is equal to the periodic direction field when  $(\phi_q^p,\phi_s^r)=0$ for $s\ne q$.
\end{corollary}

Let us say that two indices $s$ and $q$ are a {\em critical pair} if for some indices $r$ and $p$
the quantity $(\phi^r_s, \phi^p_q)=\langle\phi^r_s dx, \phi^p_q dy\rangle$ is non-zero.

\begin{corollary} 
\label{lemma: AP,k} The periodic direction field $K$ is the set $\{ \lambda \in k : \phi_s(\lambda)=\phi_q(\lambda)\} $ for all critical pairs $s$ and $q$ with $s\ne q$. 
\end{corollary}

A theorem of Calta in (\cite{calta}) shows that the collection of slopes of algebraically periodic directions of an algebraically periodic surface in standard form is a field. Note that the characterization given above gives an alternate proof that these slopes form a field.

\begin{proposition} \label{proposition:SA} Let $\lambda \in k$. Then $\lambda \in K$, iff multiplication by $\lambda$ is a self-adjoint transformation of $\Gamma^*$  with respect to $(\cdot, \cdot)$. \end{proposition}

\begin{proof} Let $\lambda \in K$ and let $\phi_q^p$ and $\phi_s^r$ be elements of $\Gamma^*$. Recall that if $q,s$ is a critical pair, then $\phi_s(\lambda)=\phi_q(\lambda)$. If $q,s$ is a critical pair, we have 
\begin{align*} (\lambda \cdot \phi_q^p, \phi_s^r) & = \langle \phi_q(\lambda) \phi_q^p dx , \phi_s^r dy\rangle\\ 
&= \langle \phi_s(\lambda) \phi_q^p dx , \phi_s^r dy\rangle  \\ 
&= \langle  \phi_q^p dx , \phi_s(\lambda) \phi_s^r dy\rangle  \\
&=(\phi_q^p, \lambda \cdot \phi_s^r)
\end{align*}

On the other hand, if $q,s$ is not a critical pair, then each of the integrals in the above equations is zero. In either case,
$(\lambda \cdot \phi_q^p , \phi_s^r)=(\phi_q^p, \lambda \cdot \phi_s^r)$. The proof of the backward implication is essentially the same. This proves the proposition. 
\end{proof} 

The construction of the maps $\phi_q^p$ began with the construction of $k$-linear maps from $\Gamma$ to $k$. Since $K\subset k$ we can also view $\Gamma$ as a $K$-vector space. We will now construct a $\QQ$ basis for $\Hom_\QQ(\Gamma,\CC)$ based on $K$-linear maps.
Let $\eta_q^p=\eta_q \circ \eta^p$ where $\eta_q : K \to \CC$ is a field embedding and $\eta^p : \Gamma \to K$ is $K$-linear.

\begin{lemma} \label{lemma:perp} 
Let $J$ be algebraically periodic and in standard form. Consider two $\QQ$-linear maps $\eta_q^p$ and $\eta_s^r$ of $\Gamma$ into $\CC$.   If $q\ne s$ then $(\eta_q^p,\eta_s^r)=0$. \end{lemma}

\begin{proof} 
We wish to show that when $q \neq s$, that is, when $\eta_q \neq \eta_s$, then $(\eta_q^p, \eta_s^r)=0$ for any maps $\eta^p$ and $\eta^r$. 
If $\eta_q \neq \eta_s$, then there exists a $\lambda \in K$ such that $\eta_q(\lambda) \neq \eta_s(\lambda)$. We have 
$$  (\lambda \cdot \eta_q^p, \eta_s^r)=  \eta_q(\lambda) (\eta_q^p, \eta_s^r)$$ 
and 
$$(\eta_q^p, \lambda \cdot \eta_s^r)=  \eta_s(\lambda) (\eta_q^p, \eta_s^r).$$ 

Since $\lambda$ acts self-adjointly, we have that $(\lambda \cdot \eta_q^p, \eta_s^r)=(\eta_q^p, \lambda \cdot \eta_s^r)$. This, in combination with the above two equations implies that either $\eta_q(\lambda)=\eta_s(\lambda)$ or $(\eta_q^p,\eta_s^r)=0$. Thus, $\eta_q^p$ and $\eta_s^r$ are perpendicular. 
\end{proof} 

\begin{theorem} \label{lemma:uniqueJ}  Let $S$ be algebraically periodic. If we restrict the $\CC$-bilinear, symmetric form $(\cdot, \cdot)$ on $Hom(\Gamma, \CC)$ to $Hom_K(\Gamma, K)$, then the restricted form $(\cdot, \cdot)_K$ determines $J$. 
\end{theorem}

\begin{proof}  Let $\psi,\phi:\Gamma \to K$ be $K$-linear so that $\psi, \phi \in Hom_K(\Gamma,K)$ and let $\sigma, \tau : K \to \CC$ be field embeddings. By Lemma~\ref{lemma:perp}, if $\sigma \neq \tau$, then $(\sigma \psi, \tau \phi)=0$. So suppose that $\sigma = \tau$ and write $$J=\sum_i \begin{bmatrix} a_i \\ c_i \end{bmatrix} \wedge \begin{bmatrix} b_i \\ d_i \end{bmatrix}.$$ We have:
\begin{align*} (\sigma \psi, \sigma \phi) &= \langle \sigma \psi dx , \sigma \phi dy\rangle  \\ 
&=\frac{1}{2}\sum_i \sigma \psi(a_i) \sigma \phi(d_i) - \sigma \psi(c_j) \sigma \phi(b_j) \\ 
&=\sigma \Big( \frac{1}{2} \sum_i \psi(a_i) \phi(d_i) -  \psi(c_j) \phi(b_j) \Big)  \\ 
&=\sigma(\psi,\phi) \notag \end{align*} 

So the inner product on $Hom(\Gamma, \CC)$ is determined by the restricted form on $Hom_K(\Gamma,K)$ and on the field embeddings of $K$, which are determined by $K$ itself. It follows from Proposition \ref{prop:innerprod} that the inner product determines $J$.
\end{proof} 

\begin{proposition} If $J\in\RR^2\wedge_\QQ\RR^2$ has rank one then $J$ is determined by its periodic direction field and its area.
\end{proposition}

\begin{proof} Put $J$ in standard form. If $J$ has rank one then then we can choose coordinates so that $\Gamma=K$. In this case $Hom_K(\Gamma,K)$ is the one dimensional vector space generated by the identity map 1. The form $(\cdot,\cdot)$ is determined by its value on this generator and we have
$(1,1)=\int_S dx\wedge dy=area(S).$
\end{proof}

\begin{cor:1.11-2} Let $S$ be completely algebraically periodic. Then the $J$-invariant of $S$ is completely determined by the periodic direction field and the area of $S$.
\end{cor:1.11-2}

\begin{proof}
This follows because if $S$ is completely algebraically periodic then $J(S)$ has rank one.
\end{proof}

We will give explicit formulae for these $J$ invariants in section 7.

If the area of $S$ is an element of the periodic direction field $K$ then the area can be changed without affecting the set of periodic directions by applying a linear transformation with coefficients in $K$ and determinant $1/area$. We can take for example the linear transformation 
$$\begin{bmatrix} 1/area &0\\0&1\end{bmatrix}.$$

Given a number field $K$ we define the {\em canonical $J$ invariant $J_K$ associated to the field $K$} to be the unique algebraically periodic $J$ invariant with area 1 which has periodic direction field $K$ and which can be written with coefficients in $K$. We denote this by $J_K$. In Section \ref{sec:formJ} we will derive formulae for $J_K$.

We say that two elements of $\RR^2\wedge_\QQ\RR^2$ are {\em isogenous} if one is obtained from the other by applying a linear transformation of the form:
$$\begin{bmatrix}\gamma &0\\0&\gamma\end{bmatrix}.$$

\begin{theorem} If $J$ is algebraically periodic and in standard form with periodic direction field $K$ and $J$ has rank $n$ then $J$ is the sum of $n$ isogenous copies of the canonical invariant $J_K$.
\end{theorem}
\begin{proof} Let $\eta^1, \ldots , \eta^n$ be an orthogonal basis for the form $(\cdot,\cdot)$ on $\Hom_K(\Gamma,K)$. Let $\gamma_1\ldots \gamma_n$ be a dual basis over $K$ for $\Gamma$. 
Any $c\in \Gamma$ can be written as $$c=\sum_\ell \eta^\ell(c)\gamma_\ell.$$

Since the $x$ and $y$ axes are algebraically periodic we can write $J$ as:
$$J(S)=\sum_i \begin{bmatrix} \alpha_i \\ 0 \end{bmatrix} \wedge \begin{bmatrix} 0 \\ \beta_i \end{bmatrix}$$

Corresponding to the decomposition of $\Gamma$ in terms of the orthogonal basis $\eta^1, \cdots, \eta^n$, there is a decomposition of $J$.
Define 

\begin{equation} \label{eq:decomp} J_{\ell,m}=\sum_i \begin{bmatrix} \eta^\ell(\alpha_i)\gamma_\ell \\ 0 \end{bmatrix} \wedge \begin{bmatrix} 0 \\ \eta^m(\beta_i)\gamma_m \end{bmatrix}. \end{equation} 

Using bilinearity to expand $J$ gives:
 $$J=\sum_{\ell,m} J_{\ell,m}.$$
 Note that $J_{\ell,\ell}$ is isogenous to a $J_{\ell,m}$ with coefficients in $K$. We wish to show that the terms $J=\sum_{\ell,m}$ with $\ell\ne m$ vanish. We will do this by looking at the corresponding forms $\langle\cdot,\cdot\rangle$. 

The value of $J$ is determined by the form $\langle\cdot,\cdot\rangle$. 
Equation~\ref{eq:decomp} shows that the horizontal and vertical directions are algebraically periodic. Recall that for $\sigma,\tau \in \Lambda^*$, 

\begin{align}\label{eq:parang} (\sigma,\tau) &= \langle \sigma dx, \tau dy \rangle 
&= \frac{1}{2} \sum_i \sigma(a_i)\tau(d_i)-\sigma(c_i) \tau(b_i)
\end{align}

where $$J=\sum_i \begin{bmatrix} a_i \\ b_i \end{bmatrix} \wedge \begin{bmatrix} c_i \\ d_i \end{bmatrix}.$$ 

In our case, Equation~\ref{eq:decomp} and Equation~\ref{eq:parang}show tha

\begin{align*}
\langle \eta^p_q dx, \eta^r_s dy\rangle&=\frac{1}{2}\sum_k \eta^p_q(\alpha_k)\eta^r_s(\beta_k).
\end{align*}
We define $\langle\cdot,\cdot\rangle_{\ell,m}$ to be the skew-symmetric bilinear form corresponding to $J_{\ell,m}$.
We calculate that 
\begin{align*}
\langle \eta^p_q dx, \eta^r_s dy\rangle_{\ell,m}&=\frac{1}{2}\sum_k \eta^p_q(\eta^\ell(\alpha_k)\gamma_\ell)
\eta^r_s(\eta^m(\beta_k)\gamma_m)\\
&=\frac{1}{2} \sum_k \eta^\ell_q(\alpha_k)\eta^p_q(\gamma_\ell)\eta^m_s(\beta_k)\eta^r_s(\gamma_m).
\end{align*}

Note that $\eta^p_q(\gamma_\ell)$ is 0 unless $p=\ell$. So if $p=\ell$, we have

 \begin{align*}
\langle \eta^\ell_q dx, \eta^m_s dy \rangle_{\ell,m}&=\frac{1}{2} \sum_k \eta^\ell_q(\alpha_k)\eta^\ell_q(\gamma_\ell)\eta^m_s(\beta_k)\eta^m_s(\gamma_m)\\
&=\frac{1}{2}\sum_k \eta^\ell_q(\alpha_k)\eta^m_s(\beta_k)\\
&=\langle \eta^\ell_q dx, \eta^m_s dy \rangle.
\end{align*}

The quantity $\langle \eta^\ell_q dx, \eta^m_s dy \rangle$ 
is 0 unless $q=s$ by Lemma \ref{lemma:perp} and $\langle \eta^\ell_q, \eta^m_q\rangle$ is 0 unless $\ell=m$ since the $\eta^1\ldots \eta^n$ are orthogonal. We conclude that $\langle\cdot,\cdot\rangle_{\ell,m}$ is zero for $\ell\ne m$. It follows that $J_{\ell,m}=0$ for $\ell\ne m$.
Thus $$J=\sum_{m} J_{m,m}.$$
It remains to show that  $J_{m,m}$ is algebraically periodic and that the periodic direction field is $K$.
This follows from Corollary~\ref{cor:k=K} since $(\eta^m_q, \eta^m_s) = \langle \eta^m_q dx, \eta^m_s dy \rangle=0$ for $q\ne s$.
\end{proof}

\begin{thm:1.6-2}
Assume that $a_j, b_j, c_j,$ and $d_j$ lie in a number field $L$. Then the lines of slope 0, 1 and $\infty$ are algebraically periodic directions if and only if  for any pair of distinct field embeddings $\sigma, \tau$ of $L$ into
$\CC$, the following equations hold: 

\[ \sum_j \sigma(c_j) \tau(b_j)- \sigma(b_j) \tau(c_j)= \sum_j \sigma(d_j) \tau(a_j)- \sigma(a_j) \tau(d_j)
\]

\[
\sum_j \sigma(c_j) \tau(d_j)- \sigma(d_j) \tau(c_j)= \sum_j
\sigma(a_j)\tau(b_j)-\sigma(b_j)\tau(a_j)=0
\]
Furthermore if the above equations hold, then the periodic direction field is the field consisting of $\lambda\in L$  for which $\sigma(\lambda)=\tau(\lambda)$ whenever $\sum_j \sigma(c_j) \tau(b_j)- \sigma(b_j) \tau(c_j)\ne0$.
\end{thm:1.6-2}

%\end{thm:1.9-2}

\begin{proof} This follows from Theorem \ref{alg per crit} and Corollary \ref{lemma: AP,k} and the observations that the field embeddings of $L$ give a basis for $Hom_\QQ(L,\CC)$.
\end{proof}

\end{section}

 \begin{section}{Formulae for J}
\label{sec:formJ}
 
 As we have seen the periodic direction field of a surface $S$ depends only on $J(S)$. In this section we start with a field $K$ and find an element $J\in \RR^2\wedge_\QQ\RR^2$ which can be expressed with coefficients in $K$ and which has periodic direction field $K$. We will use this construction in section \ref{sec:construction} where we will construct surfaces  which realize these invariants.

In this section we will deal with $J$ which are written as 
\begin{equation}
\label{eq:diagform}
J=\sum_{i=1}^{n} \begin{bmatrix}a_i \\ 0
\end{bmatrix} \wedge \begin{bmatrix}0 \\ b_i \end{bmatrix}.
\end{equation}
This representation implies that the $x$ and $y$ axes are algebraically periodic.
Conversely if the $x$ and $y$ axes are algebraically periodic then $J$ can be written in this form.
Before we begin we record a fact that will be useful in Section~\ref{sec:decomp}.

\begin{proposition}
\label{prop:n ge d}
If $J$ is algebraically periodic with periodic direction field $K$ of degree $d$ and can be written as above where $n$ is the number of terms in the expression then $n\ge d$.
\end{proposition}
\begin{proof} 
By using a coordinate change that preserves the axes we may assume that $J$ is in standard form. Assume that $J$ can be written as above with $n<d$. Since $\Gamma$ is a $K$ vector space it has $\QQ$-dimension at least $d$. Thus we can choose a non-zero $\QQ$-linear map $\tau:\Gamma\to\QQ$ so that $\tau(b_i)=0$ for $i=1\ldots n$.
Now for any $\sigma:\Gamma\to\QQ$ we have 
$$(\sigma,\tau)=\frac{1}{2}\sum_{i=1}^n \sigma(a_i)\tau(b_i)=0.$$
But this contradicts the non-degeneracy of $(\cdot,\cdot)$ proved in Corollary~\ref{cor:symform}.
\end{proof}

We record a useful criterion for an element $J\in \RR^2\wedge_\QQ\RR^2$ to have rank one. This is a special case of Theorem~\ref{thm:1.6-1}.

\begin{lemma}
\label{lem:criterion}
Say that $J$ can be written as in Equation~\ref{eq:diagform}
with $a_i$ and $b_i$ in the field $K$. Let $\sigma_j$ be the field 
embeddings of $K$ into $\CC$. The periodic direction field of $J$ will be $K$ exactly when:
$$
\sum_{i=1}^n\sigma(a_i)\tau(b_i)=0
$$
for any pair of distinct field embeddings $\sigma$ and $\tau$.
\end{lemma}

\begin{proof} This follows from Theorem~\ref{thm:1.6-1}.
\end{proof}

The next three theorems give explicit expressions for $J$ when $J$ has rank one  with periodic direction field $K$. The first is in terms of an integral matrix with an eigenvalue that generates $K$. The second is in terms of a minimal polynomial of $\lambda$ where $\lambda$ generates $K$. The third is in terms of the trace of $K$.

\begin{theorem} \label{thm:integral} Let $A$ be an $n \times n$ integral matrix with eigenvalue
$\lambda$ and corresponding eigenvector $(u_i)$. Then $\lambda$ is
an eigenvalue for $A^{t}$ as well. Let $(v_i)$ be a corresponding
eigenvector. Then if 
$$J=\sum_{i=1}^{n} \begin{bmatrix}u_i \\ 0
\end{bmatrix} \wedge \begin{bmatrix}0 \\ v_i \end{bmatrix},$$
 $J$ is algebraically periodic with periodic direction field $K={Q}(\lambda)$.
\end{theorem}

\begin{proof} We show that the periodic direction field of $J$ is $K$.
We begin by showing that
\begin{equation}
\sum_{i=1}^{n} \begin{bmatrix} \lambda u_i \\ 0 \end{bmatrix} \wedge
\begin{bmatrix} 0 \\ v_i \end{bmatrix} 
=
 \sum_{i=1}^{n} \begin{bmatrix}u_i \\ 0 \end{bmatrix} \wedge
\begin{bmatrix} 0 \\ \lambda v_i. \end{bmatrix}.
\end{equation}

Note that $\lambda u_i=\sum_{j=1}^{n}a_{ij}u_j$ and $\lambda v_i = \sum_{j=1}^{n}a_{ji}
v_j$. This implies that
\begin{align*}
\sum_{i=1}^{n} \begin{bmatrix} \lambda u_i \\ 0 \end{bmatrix} \wedge
\begin{bmatrix} 0 \\ v_i \end{bmatrix} 
&= 
\sum_{i=1}^{n} \sum_{j=1}^{n} a_{ij} \begin{bmatrix}u_j \\ 0 \end{bmatrix} \wedge \begin{bmatrix} 0 \\ v_i \end{bmatrix}\\
& = \sum_{i=1}^{n} \sum_{j=1}^{n} a_{ji} \begin{bmatrix}u_i \\ 0 \end{bmatrix} \wedge \begin{bmatrix}0 \\ v_j \end{bmatrix}\\
&=
 \sum_{i=1}^{n} \begin{bmatrix}u_i \\ 0 \end{bmatrix} \wedge
\begin{bmatrix} 0 \\ \lambda v_i \end{bmatrix}.
\end{align*}

Now consider two distinct field embeddings $\sigma$ and $\tau$. Since $\lambda$ generates $K$ it must be the case that $\sigma(\lambda)\ne\tau(\lambda)$.
\begin{align*}
\sigma(\lambda)\sum_{i=1}^n\sigma(u_i)\tau(v_i)&=
\sum_{i=1}^n\sigma(\lambda u_i)\tau(v_i)\\
&=\sum_{i=1}^n\sigma(u_i)\tau(\lambda v_i)\\
&=\tau(\lambda)\sum_{i=1}^n\sigma(u_i)\tau(v_i)
\end{align*}
Thus $\sum_{i=1}^n\sigma(u_i)\tau(v_i)=0.$ By Lemma \ref{lem:criterion}, this proves the theorem. 
\end{proof}

\begin{thm:1.12-2} If $\alpha_1\ldots\alpha_n$ and $\beta_1\ldots\beta_n$ are  dual bases for the trace form on a field $K$ then 
$$J=\sum_j \begin{bmatrix}\alpha_j\\0\end{bmatrix}\wedge \begin{bmatrix}0\\ \beta_j\end{bmatrix}$$ is algebraically periodic with periodic direction field $K$. Conversely if $J$ can be written as above and with entries in $K$ and $J$ is algebraically periodic with periodic direction field $K$ and area 1 then $\alpha_1\ldots\alpha_n$ and $\beta_1\ldots\beta_n$ are a dual basis for the trace form on $K$.
\end{thm:1.12-2}

\begin{proof} Let $\phi_1\ldots\phi_n$ be field embeddings of $K$ into $\CC$.
Define $n\times n$ matrices $A$ and $B$ by setting
$a_{ij}=\phi_i(\alpha_j)$ and $b_{ij}=\phi_i(\beta_j)$. The fact that the $\alpha$ and $\beta$ are dual bases for the trace translates into the fact that $A^TB=I$:

\begin{align*}
[A^TB]_{ik}&=\sum_j [A^T]_{ij}[B]_{jk}=\sum_j \phi_j(\alpha_i)\phi_j(\beta_k)\\
&=\sum_j \phi_j(\alpha_i \beta_k)\\
&=tr(\alpha_i\beta_k)\\
&=\delta_{ik}.
\end{align*}

Consider the inner product corresponding to $J$. 
We now want to calculate $(\phi_i,\phi_k)$ and show that $(\phi_i,\phi_k)=\delta_{ik}$. This will establish that $J$ is the element corresponding to the field $K$.

$$(\phi_i,\phi_k)=\langle\phi_i dx,\phi_k dy\rangle=\sum_j\phi_i(\alpha_j)\phi_k(\beta_j).$$

$$\sum_j\phi_i(\alpha_j)\phi_k(\beta_j)=\sum_j [A]_{ij}[B^T]_{jk}=[AB^T]_{ik}$$

Since $B$ is a right inverse to $A^T$, $B$ is also a left inverse to $A^T$ so $BA^T=I$. Taking transposes of both sides gives $AB^T=I$. Thus $(\phi_i,\phi_k)=\delta_{ik}$. According to Lemma~\ref{lem:criterion} $J$ has rank one with periodic direction field $K$. We note that this $J$ will have area 1.
\end{proof}

The next result gives an explicit way to write $J=J_K$ in terms of the field $K$.

\begin{theorem} \label{thm:J(k)} Let $\lambda$ be an algebraic number with minimal polynomial $$p(x)=x^n + a_{n-1}x^{n-1} + \cdots + a_1x + a_0$$ and let $k=\QQ(\lambda)$. Let
$$q(x)=\frac{p(x)}{x-\lambda}=b_{n-1}x^{n-1} + b_{n-2}x^{n-2} + \cdots + b_1x + b_0.$$
Let $\beta_j= \frac{b_j}{p'(\lambda)}$.
Then
$$J_K= \sum_{j=0}^{n-1} \begin{bmatrix} \lambda^j \\ 0 \end{bmatrix} \wedge \begin{bmatrix} 0 \\ \beta_j \end{bmatrix}.$$

\end{theorem}

\begin{proof}
This follows since $\lbrace \lambda^j \rbrace_{j=0}^{n-1}$ and $\lbrace \frac{b_j}{p'(\lambda)} \rbrace_{j=0}^{n-1}$ are dual bases for the trace pairing, $tr: k \times k \to \QQ$. (See \cite{lang}, Proposition $1$, pg. 213.) 
\end{proof}

\end{section}

\begin{section}{$\Iso(J)$ in the algebraically periodic case}

We now calculate $\Iso(J)$ in the algebraically periodic case. Note that we are not assuming that $J$ has rank one. The fact that $\Iso(J)$ acts transitively on the set of algebraically periodic directions was proved in (\cite{calta}). 

\begin{theorem}  If $J$ is algebraically periodic and in standard form with periodic direction field $K$ then $ \Iso(J)=SL(2,K)$.
\label{thm:isoJcalc}
\end{theorem}

\begin{proof}

Let $$M=\begin{bmatrix} \alpha & \beta \\ \gamma & \delta \\ \end{bmatrix}$$
be an element of $SL(2,k)$. We compute the action of this linear transformation on the space $\Lambda^*$.

The symbol $\phi^k_j dx \in \Lambda^*$ represents the homomorphism which takes a vector $\begin{bmatrix}x\\y\end{bmatrix}$ in $\Lambda$ to the complex number $\phi^k_j(x)$. If we precompose with the matrix $M$ we get the homomorphism which takes the vector $\begin{bmatrix}x\\y\end{bmatrix}$ to $\begin{bmatrix}\alpha x+\beta y\\ \gamma x +\delta y\end{bmatrix}$
and then to $\phi^k_j(\alpha x+\beta y)$.
Now $\phi^k_j(\alpha x+\beta y)=\phi_j(\alpha)\phi_j^k(x)+\phi_j(\beta)\phi_j^k(y)$. 
If we express this form in $dx$ and $dy$ notation it becomes:
$\phi_j(\alpha)\phi_j^k dx+\phi_j(\beta)\phi_j^k dy$.

A similar calculation shows that $\phi^k_j dy$ is taken to 
$\phi_j(\gamma)\phi_j^k dx+\phi_j(\delta)\phi_j^k dy$.

Corollary~\ref{thm:HSG} implies that $J$ will be preserved if and only if integrals of products of holonomic 1-forms are preserved. It suffices to check this on our basis of 1-forms.

This will occur when 
$$\langle\phi^p_q dx, \phi^r_s dy\rangle=\langle \big(\phi_q(\alpha)\phi^p_q dx+\phi_q(\beta)\phi^p_q dy\big), \big(\phi_s(\gamma)\phi^r_s dx+\phi_s(\delta)\phi^r_s dy\big)\rangle$$

Expanding the right hand side we have:

\begin{multline} 
\label{eq:phipq} 
\langle \phi_q(\alpha)\phi^p_q dx,\phi_s(\gamma)\phi^r_s dx\rangle
+\langle \phi_q(\alpha)\phi^p_q dx, \phi_s(\delta)\phi^r_s dy\rangle\\
+\langle \phi_q(\beta)\phi^p_q dy, \phi_s(\gamma)\phi^r_s dx\rangle
+\langle \phi_q(\beta)\phi^p_q dy, \phi_s(\delta)\phi^r_s dy\rangle
\end{multline}

Theorem \ref{alg per crit} shows that the $\langle dx, dx\rangle $ terms and the $\langle dy, dy\rangle$ terms vanish so expression~\ref{eq:phipq} becomes:

\begin{equation}
\begin{split}
&\langle \phi_q(\alpha)\phi^p_q dx, \phi_s(\delta)\phi^r_s dy\rangle
+\langle \phi_q(\beta)\phi^p_q dy, \phi_s(\gamma)\phi^r_s dx\rangle\\
&=\phi_q(\alpha) \phi_s(\delta)\langle \phi^p_q dx, \phi^r_s dy\rangle
+\phi_q(\beta)\phi_s(\gamma)\langle \phi^p_q dy, \phi^r_s dx\rangle\\
&=\big(\phi_q(\alpha) \phi_s(\delta)-\phi_q(\beta)\phi_s(\gamma)\big)
\langle \phi^p_q dx, \phi^r_s dy\rangle
\end{split}
\end{equation}

Summarizing, the integral of the wedge product is invariant when
\begin{equation}\langle \phi^p_q dx, \phi^r_s dx\rangle=\big(\phi_q(\alpha) \phi_s(\delta)-\phi_q(\beta)\phi_s(\delta)\big)
\langle \phi^p_q dx, \phi^r_s dy\rangle
\end{equation}
This equation will hold when 
\begin{equation}
\phi_q(\alpha) \phi_s(\delta)-\phi_q(\beta)\phi_s(\gamma)=1
\end{equation}
or when $\langle\phi^p_q dx, \phi^r_s dy\rangle=0$.
Thus we have invariance of $J$ if and only if for any critical pair of field embeddings (not necessarily distinct) we have 
$$\phi_q(\alpha) \phi_s(\delta)-\phi_q(\beta)\phi_s(\gamma)=1$$

Now suppose that $M \in \Iso(J)$. We will show that $M=\lambda M'$ for some $M' \in GL(2,K)$ and $\lambda \in k$ such that 
$$
M' = \begin{bmatrix} a & b \\ c & d \\ \end{bmatrix}
.$$
We know that $M \in SL(2,k)$ and that $M$ preserves the algebraically periodic directions. Thus $M$ takes the triple of directions $(0,1,\infty)$ to another triple of directions $(s_1, s_2, s_3)$ where each $s_i \in K$. There exists a matrix $M' \in GL(2,K)$ which takes $(0,1,\infty)$ to $(s_1, s_2, s_3)$.
Therefore $(M')^{-1} M$ takes the triple $(0,1,\infty)$ to itself. This implies that $(M')^{-1} M = \lambda I$ where $\lambda \in k$.  

We consider equation \ref{deteqn} for $M$ in this special form.
$$ \phi_q(\lambda a) \phi_s(\lambda d)-\phi_q(\lambda b)\phi_s(\lambda c)=1$$
for any critical pair of field embeddings.

This equation is equivalent to:

\begin{equation}  \label{deteqn} \phi_q(\lambda)\phi_s(\lambda)\big(\phi_q( a) \phi_s( d)-\phi_q( b)\phi_s(c)\big)=1 \end{equation} 

Note that Corollary~\ref{lemma: AP,k} implies that $\phi_q$ and $\phi_s$ agree on elements of $K$ so we get

\begin{equation}
\begin{split}
\phi_q( a) \phi_s( d)-\phi_q( b)\phi_s(c)&=\phi_q( a) \phi_q( d)-\phi_q( b)\phi_q(c)\\
&=\phi_q(ad-bc)=\phi_s(ad-bc)
\end{split}
\end{equation}

Now, if $\lambda \in K$, then $\phi_q(\lambda)=\phi_s(\lambda)$, and so equation~\ref{deteqn} holds. We will now assume equation~\ref{deteqn} holds and show that $\lambda \in K$. 

Write $\theta$ for $(ad-bc)^{-1}$. Equation \ref{deteqn} holds if $\phi_q(\lambda)\phi_s(\lambda)=\phi_q(\theta)=\phi_s(\theta)$.
If we take $\phi_q=\phi_s=id$ we see that $\lambda^2=\theta$. 
It follows that $\phi_q(\lambda)\phi_s(\lambda)=\phi_q(\lambda^2)=\phi_q(\theta)$. Thus, $\phi_q(\lambda)=\pm\sqrt{\phi_q(\theta)}$. Similarly, $\phi_s(\lambda)=\pm\sqrt{\phi_s(\theta)}$.

If $\phi_q(\lambda)=\phi_s(\lambda)$ then $\phi_q(\lambda)\phi_s(\lambda)=\phi_q(\theta)$. On the other hand, if $\phi_q(\lambda)\ne\phi_s(\lambda)$ then $\phi_q(\lambda)\phi_s(\lambda)=-\phi_q(\theta)$.
Thus Equation~\ref{deteqn} holds exactly when $\phi_q(\lambda)=\phi_s(\lambda)$. Since this equation holds for any critical pair then we conclude that $\lambda\in K$ and $M=\lambda M'\in SL(2,K)$.
\end{proof}

\begin{corollary} \label{corollary:trace} Let $J$ be algebraically periodic. Each element of $\Iso(J)$ has trace in the periodic direction field of $J$. \end{corollary}

 \begin{corollary}\label{corollary: cont}  Let $J$ be algebraically periodic. Let $K$ be the periodic direction field, $k$ the multiplicative field, $L$ the holonomy field and $k'$ the trace field of $\Iso(J)$. Then $k' \subset K  \subset k \subset L$.
 \end{corollary}

  \begin{proof} The first inclusion follows from the preceding corollary, the second from Proposition~\ref{prop:apk}, and the last from Proposition~\ref{prop: k,K}.
  \end{proof}
  
 \end{section}
 
 \begin{section} {Surfaces with affine automorphisms and the homological affine group}
\label{section:reptheory}

In section $6$ we found formulae which can be used to detect the property of algebraic periodicity. In this section we  show that the existence of certain affine automorphisms of $S$ imply that $S$ is algebraically periodic. 

\begin{thm:1.2-2} If the Veech group of $S$ contains a  pseudo-Anosov element $A$, then the periodic direction field $K$ of $S$ is equal to the trace field and $S$ is completely algebraically periodic. \end{thm:1.2-2} 

\begin{proof} By replacing $A$ by its square if necessary, we may assume that $A$ preserves the orientations on the stable and unstable foliations. 

The following argument is from McMullen (see \cite{ctm:dynamics}).
We can choose a basis of $EH^*$ consisting of forms $\sigma_j dx$ and $\sigma_j dy$ where $\sigma_j$ runs through the field embeddings of $k'$. 
Let $M\in\slr$ be the matrix of the pseudo-Anosov element. Since $M$ preserves $J(S)$ we have $\langle M\alpha, M\beta\rangle=\langle\alpha,\beta\rangle$ for any $\alpha$ and $\beta$ in $EH^*$. It follows that 

\begin{equation*}
\begin{split}
\langle M\alpha,\beta\rangle&=\langle M\alpha,MM^{-1}\beta\rangle\\
&=\langle\alpha,M^{-1}\beta\rangle,\\
\langle M^{-1}\alpha,\beta\rangle&=\langle M^{-1}\alpha,M^{-1}M\beta\rangle\\
&=\langle\alpha,M\beta\rangle,\\
\big\langle (M+M^{-1})\alpha,\beta\big\rangle&=\big\langle\alpha,(M+M^{-1})\beta\big\rangle.
\end{split}
\end{equation*}

Since $\det M=1$ it follows that $M+M^{-1}=tr(M)I$.

Now applying the transformation $M+M^{-1}$ to our basis elements gives 

\begin{align*}
\sigma_j(tr M)\langle\sigma_j dx,\sigma_k dy\rangle
&=\langle(M+M^{-1})\sigma_j dx,\sigma_k dy\rangle\\
&=\langle\sigma_j dx,(M+M^{-1})\sigma_k dy\rangle\\
&=\sigma_k(tr M)\langle\sigma_j dx,\sigma_k dy\rangle
\end{align*}

Since $tr M$ generates $k'$ we have $\sigma_j(tr M)\ne \sigma_k(tr M)$ unless $j=k$.
We conclude that $\langle\sigma_j dx,\sigma_k dy\rangle=0$ when $j\ne k$.
The same argument also gives  $\langle\sigma_j dx,\sigma_k dx\rangle=0$ and $\langle\sigma_j dy,\sigma_k dy\rangle=0$ when $j\ne k$. Theorem~\ref{thm:1.6-1} shows that $S$ is algebraically periodic and that the periodic direction field of $S$ is the trace field.

According to (\cite{kenyon-smillie}) the trace field is equal to the holonomy field. Since the holonomy field is equal to the periodic direction field $S$ is completely algebraically periodic according to Corollary~\ref{cor:capeq}.
\end{proof} 

\begin{cor:1.3-2} If a flat surface $S$ is obtained from the Thurston-Veech construction then it is completely algebraically periodic.
\end{cor:1.3-2}
\begin{proof}
In this case the Veech group contains a pseudo-Anosov automorphism.
\end{proof}

What we have actually shown is the following:

\begin{theorem} If there is an affine automorphism of $S$ with trace $\theta$ and $\theta$ generates the holonomy field then $S$ is completely algebraically periodic.
\end{theorem}

\begin{thm:1.4-2} If $S$ is a flat surface obtained from a rational-angled triangle via the Zemlyakov-Katok construction then $S$ is completely algebraically periodic.
\end{thm:1.4-2}

\begin{proof} Let $\Delta$ be the triangle with angles $(p_1/q)\pi$, $(p_2/q)\pi$ and $(p_3/q)\pi$ where $p_1$, $p_2$ and $p_3$ have no common factor. Assume that the segment between $(0,0)$ and $(1,0)$ is one edge of the triangle.
Let $\theta=e^{2\pi i/q}$. Let $\QQ(\theta)$ be the cyclotomic field generated by $\theta$. 
If we identify the plane $\RR^2$ with $\CC$ then we can see that the reflections in the sides of the triangle are contained in the group generated by rotations through the angle of $2\pi/q$, complex conjugation and translation by elements in $\QQ(\theta)$. It follows that all triangles obtained by repeatedly unfolding $\Delta$ have coordinates in $\QQ(\theta)$. If we return to viewing $\CC$ as $\RR^2$ then the vertices of these triangles are vectors and the coordinates of these vectors lie in the collection of real and imaginary parts of elements of $\QQ(\theta)$. This is just the real field $\QQ(\theta+\theta^{-1})$. In particular the holonomy of $S$ is contained in $\QQ(\theta+\theta^{-1})$ by Proposition~\ref{prop:hol criterion}.
 The matrix that corresponds to rotation by $2\pi/q$ induces an automorphism of the surface $S$ which has trace $2\cos(2\pi/q)=\theta+\theta^{-1}$. 
In this case the holonomy field is contained in the trace field.  And since the trace field is contained in the holonomy field, it follows that $k'=K=L$ and so  $S$ is completely algebraically periodic.
\end{proof}

\begin{theorem} If $\Iso(J)$ contains a  rational matrix which is diagonalizable over the rationals with eigenvalues other than 1 and -1 then the eigendirections are algebraically periodic.
\label{thm:eigen}
\end{theorem}

\begin{proof} Let $M \in \Iso(J)$. We can choose coordinates so that $M$ is diagonal with diagonal entries $\lambda$ and $\lambda^{-1}$. We can write $J$ as:
$$J= \sum_j \begin{bmatrix} a_j \\ 0 \end{bmatrix} \wedge \begin{bmatrix} 0 \\ b_j \end{bmatrix}+\sum_k \begin{bmatrix} c_k \\ 0 \end{bmatrix} \wedge \begin{bmatrix} d_k \\ 0\end{bmatrix}
+\sum_\ell \begin{bmatrix} 0\\ e_\ell \end{bmatrix} \wedge \begin{bmatrix} 0 \\ f_\ell\end{bmatrix}.
$$

We get 

$$MJ= \sum_j \begin{bmatrix} \lambda a_j \\ 0 \end{bmatrix} \wedge \begin{bmatrix} 0 \\ \lambda^{-1}b_j \end{bmatrix}+\sum_k \begin{bmatrix} \lambda c_k \\ 0 \end{bmatrix} \wedge \begin{bmatrix} \lambda d_k \\ 0\end{bmatrix}
+\sum_\ell \begin{bmatrix} 0\\ \lambda^{-1} e_\ell \end{bmatrix} \wedge \begin{bmatrix} 0 \\ \lambda^{-1}f_\ell\end{bmatrix}.
$$

Using $\QQ$ linearity we have:

$$MJ= \sum_j \begin{bmatrix} a_j \\ 0 \end{bmatrix} \wedge \begin{bmatrix} 0 \\ b_j \end{bmatrix}+\lambda^2\sum_k \begin{bmatrix} c_k \\ 0 \end{bmatrix} \wedge \begin{bmatrix} d_k \\ 0\end{bmatrix}
+\lambda^{-2}\sum_\ell \begin{bmatrix} 0\\ e_\ell \end{bmatrix} \wedge \begin{bmatrix} 0 \\ f_\ell\end{bmatrix}.
$$

Comparing the formulae we see that 
$$\begin{bmatrix} c_k \\ 0 \end{bmatrix} \wedge \begin{bmatrix} d_k \\ 0\end{bmatrix}=
\begin{bmatrix} 0\\ e_\ell \end{bmatrix} \wedge \begin{bmatrix} 0 \\ f_\ell\end{bmatrix}=0.$$
\end{proof}

The following result shows that the property of being algebraically periodic is characterized by the group $\Iso(J)$.

\begin{theorem} If $\Iso(J)$ contains a subgroup isomorphic to $SL(2,\QQ)$ as an abstract group then $J$ is algebraically periodic.
\label{thm: absgroup}
\end{theorem}

\begin{proof} The group $SL(2,\QQ)$ contains elements

$$\alpha=\begin{bmatrix} 2 &0\\0&2^{-1}\end{bmatrix}\qquad
\beta=\begin{bmatrix} 1 &1\\0&1\end{bmatrix}$$
where $\beta$ has infinite order and which satisfy $\alpha\beta\alpha^{-1}=\beta^4$.
By assumption $\Iso(J(S))$ contains elements $\alpha'$ and $\beta'$ which
satisfy this equation and for which $\beta'$ has infinite order. It follows that $\beta'$ must be a parabolic element (cf \cite{o'meara} p. 67) and the axis of $\beta'$ must be an eigenspace of $\alpha'$ with eigenvalue 2. According to Theorem \ref{thm:eigen} the axis of $\beta'$ must be an algebraically periodic direction. 

To show the existence of three algebraically periodic directions for $S$ we note that parabolic elements have a common axis if and only if they commute. Thus find three parabolic elements $\beta_j$, $j=1,2,3$ and corresponding $\alpha_j$ so that the axes of the $\beta_j$ are distinct. This implies that no pair of the $\beta$'s commute. Consider the images of these elements $\beta'_j$. These elements are parabolic and do not commute. It follows that there axes are distinct. By the above argument these axes are algebraically periodic directions.
\end{proof}

\begin{thm:1.7-2} If $\Iso(J(S))$ is isomorphic to $SL(2,K)$ as an abstract group for some field $K$ then $S$ is algebraically periodic and $K$ is the periodic direction field.
\end{thm:1.7-2}

\begin{proof} If $\Iso(J(S))$ is isomorphic to $SL(2,K)$ as an abstract group it certainly contains a subgroup isomorphic to $SL(2,\QQ)$. So by Theorem~\ref{thm: absgroup} $S$ is algebraically periodic. By Theorem \ref{thm:isoJcalc}
 $\Iso(J(S))$ is isomorphic to $SL(2,K)$ where $K$ is the periodic direction field. If $\Iso(J(S))$ were isomorphic to $SL(2,K')$ for some field $K'$ other than the direction field then we would have an isomorphism between $SL(2,K)$ and $SL(2,K')$. According to O'Meara (\cite{o'meara}) such an isomorphism can only occur when $K=K'$.
\end{proof}

Recall that $\HAG$ is the homological affine group and is defined in the introduction to this paper.
If we let $V(S)$ be the Veech group of $S$, then  $$V(S) \subset \HAG \subset \Iso(J(S)).$$

\begin{thm:1.8-2} If $S$ is completely algebraically periodic, in particular, if the Veech group of $S$ contains a hyperbolic element, then $\HAG=\Iso(J)=SL(2,K)$. \end{thm:1.8-2}

\begin{proof} Note that $\Iso(J)$ preserves the algebraically periodic directions of $S$, while $\HAG$ preserves the homological directions. When $S$ is completely algebraically periodic, (i.e., when $L=K$), the set of homological directions and the set of algebraically periodic directions are the same. Thus $\Iso(J)=\HAG$. \end{proof}

The notions of the homological affine group and $\Iso(J)$ are closely related but not identical as we show in Theorem \ref{eh ne ah}.
\end{section}

\begin{section}{Geometric Decompositions and Total Reality}
\label{sec:decomp}

We will assume throughout this section that $S$ is algebraically periodic and in standard form. If $\Lambda$ is the essential holonomy of $S$ then, as in Section \ref{sec:pdf}, we can write $\Lambda^* = \Gamma^* \oplus \Gamma^*$ where $\Gamma^*$ is the collection of $\QQ$-linear maps from $\Gamma$ to $\CC$. As was shown in Section \ref{sec:pdf} there is a  $\CC$-bilinear symmetric form on $\Gamma^*=Hom(\Gamma,\CC)$ given by $(\sigma,\tau)=\langle \sigma dx, \tau dy\rangle$.

\begin{theorem} \label{theorem:posdef} If the restriction of $(\cdot, \cdot)$ to $Hom(\Gamma,\RR)$ is positive definite, then $K$ is totally real. \end{theorem} 

\begin{proof} Fix a $K$-linear map $\psi : \Gamma \to K$ such that $(\psi,\psi) \neq 0$. Let $\lbrace \sigma_i \rbrace_{i=1}^r$ be the set of real embeddings of $K$ into $\CC$ and $\lbrace \tau_i \rbrace_{i=r+1}^n$ the subset of complex embeddings which contains exactly one of each pair of complex conjugate embeddings $\tau$ and $\bar{\tau}$. Thus, for each $1 \leq i \leq r$, we obtain one $\QQ$-linear homomorphism of $\Gamma$ into $\CC$, $\eta_i = \sigma_i \psi$. 
And for each $r+1 \leq i \leq n$, we obtain two $\QQ$-linear homomorphisms of $\Gamma$ into $\CC$, $\eta_i=\tau_i \psi$ and $\bar{\eta}=\bar{\tau}\psi$, corresponding to an embedding and its conjugate. 

Fix an $i>r$ and let $\tau=\tau_i$. Since $\tau \psi$ and $\bar{\tau} \psi$ are distinct $\QQ$-linear homomorphisms, Lemma \ref{lemma:perp} shows that $(\tau \psi, \bar{\tau} \psi)=0$. 
We have that $\tau \psi = \Re (\tau \psi) + i \Im(\tau \psi)$.  Then using $\CC$-bilinearity of $(\cdot, \cdot)$,  

\begin{align} 
(\tau \psi, \bar{\tau} \psi) &= (\Re(\tau \psi) + i \Im(\tau \psi), \Re(\tau \psi) - i \Im(\tau \psi)) \notag \\ 
&=(\Re(\tau \psi),\Re(\tau \psi)) + (\Im(\tau \psi) , \Im(\tau \psi)) \notag \\ 
&=0 \notag
\end{align}

Either both of the terms $(\Re(\tau \psi),\Re(\tau \psi))$ and $(\Im(\tau \psi)$ , $\Im(\tau \psi))$ are zero, or one is positive and the other negative. 
This contradicts the fact that $(\cdot, \cdot)$ is positive definite. Thus, there can be no complex embeddings of $K$ into $\CC$ and so $K$ is totally real. 

\end{proof}

 \begin{thm:1.14-2} Suppose that $S$ is algebraically periodic and that $S$ has a decomposition into squares so that the directions of the sides are algebraically periodic directions. Then the periodic direction field is totally real. \end{thm:1.14-2} 

\begin{proof} We may assume that $S$ is in standard form. We have not assumed that the decomposition into squares is an edge to edge decomposition. We can make it an edge to edge decomposition by adding vertices to the sides of the squares. Adding a vertex to a horizontal side changes $J$ by terms of the form $\begin{bmatrix} a \\0 \end{bmatrix} \wedge \begin{bmatrix} b \\ 0 \end{bmatrix}$.
If $a_j$ are the side lengths of the squares that decompose $S$ then we can write $J$ as

$$J= \sum_j \begin{bmatrix} a_j \\ 0 \end{bmatrix} \wedge \begin{bmatrix} 0 \\ a_j \end{bmatrix}+\sum_k \begin{bmatrix} c_k \\ 0 \end{bmatrix} \wedge \begin{bmatrix} d_k \\ 0\end{bmatrix}
+\sum_\ell \begin{bmatrix} 0\\ e_\ell \end{bmatrix} \wedge \begin{bmatrix} 0 \\ f_\ell\end{bmatrix}.
$$
Since the horizontal and vertical directions are assumed to be algebraically periodic the second and third terms are zero. So we have:
$$J= \sum_j \begin{bmatrix} a_j \\ 0 \end{bmatrix} \wedge \begin{bmatrix} 0 \\ a_j \end{bmatrix}.$$
Let $\rho:\Gamma \to \RR$ be $\QQ$-linear. Then we have that 

\begin{align*}
(\rho,\rho)&= \sum_i \rho(a_i) \rho(a_i) \\ 
&=\sum_i \rho(a_i)^2 \geq 0
\end{align*}  
Since this form is non-degenerate showing that it is non-negative shows that it is positive definite.
\end{proof}

 \begin{theorem} \label{theorem: totreal} Suppose that $S$ is algebraically periodic and that $S$ has a decomposition into rectangles so that the directions of the sides are algebraically periodic directions. Then the number of rectangles is at least the degree of the periodic direction field. \end{theorem} 
 
 \begin{proof} The proof of the previous theorem shows that the decomposition of $S$ into $n$ rectangles gives us a way to write $J$ with $n$ terms. We then apply Proposition~\ref{prop:n ge d}.
 \end{proof}

Genus two surfaces which are sums of isogenous tori play a special role (see Theorem $1.8$ of \cite{ctm:dynamics}). 

\begin{theorem} 
\label{thm:isotori} If $S$ has a decomposition into isogenous tori then $S$ is algebraically periodic and the periodic direction field is totally real.
\end{theorem}

\begin{proof} Since $J$ is additive for connected sums we can write $J(S)$ as the sum of  the $J$ invariants of tori. By choosing coordinates appropriately we may assume that the holonomy of each torus is a multiple of that of the standard integral lattice. The $J$ invariant for this torus will be 
$$ \begin{bmatrix} 1\\ 0 \end{bmatrix} \wedge \begin{bmatrix} 0 \\ 1\end{bmatrix}.$$

The surface is algebraically periodic because every direction with rational slope (including those with slopes $0, 1$ and $\infty$) is an algebraically periodic direction. So we may write 

$$J= \sum_j \begin{bmatrix} a_j \\ 0 \end{bmatrix} \wedge \begin{bmatrix} 0 \\ a_j \end{bmatrix}.$$

The result follows from the proof of Theorem~\ref{theorem:squares}. 
\end{proof}

\begin{cor:1.13-2} If $S$ has a parabolic automorphism and  $S$ is algebraically periodic then the periodic direction field $K$ is totally real. \end{cor:1.13-2}

This theorem does not assert that the existence of a parabolic automorphism implies that  $S$ is algebraically periodic. A parabolic automorphism gives rise to one algebraically periodic direction. In general we would need to have two other algebraically periodic directions to determine that $S$ was algebraically periodic but in this case the existence of a second algebraically periodic direction suffices. If there are two parabolic automorphisms with distinct axes then the total reality of $K$ follows from a result of (\cite{hl}).

\begin{proof} If $S$ has a parabolic element, then there is a complete cylinder decomposition of $S$ in the direction of the parabolic. Choose a second algebraically periodic direction and choose coordinates so that these directions are the $x$ and $y$ axes. By a coordinate change which preserves the axes we may assume that the height and width of one cylinder are both 1. For a parabolic transformation the moduli of all cylinders are rationally related so the moduli of all cylinders are rational. It follows that $S$ can be decomposed into squares. 
The result now follows from Theorem~\ref{theorem:squares}. 
\end{proof}  

The positive definiteness of the restriction of $(\cdot,\cdot)$ to $Hom_\QQ(\Gamma,\RR)$ gives non-trivial conditions for the existence of a parabolic direction or square decomposition. In (\cite{kenyon-smillie}) the positive definiteness of this inner product is used as a test for the existence of a parabolic cylinder decomposition. In that case the field is the maximal real subfield of the cyclotomic field which is totally real.

\begin{theorem} If $S$ is algebraically periodic and $S$ has a decomposition into cylinders then the number of cylinders is greater than or equal to the degree of the periodic direction field.
\end{theorem}

Note that this theorem applies to decompositions of $S$ into cylinders of arbitrary moduli.

\begin{proof} As in the previous corollary the decomposition of $S$ into cylinders gives an expression for $J$. The number of cylinders is the number of terms in this expression.
The result now follows from Proposition~\ref{prop:n ge d}.
\end{proof}

\end{section}

\begin{section}{Construction of Examples}
\label{sec:examples}
\label{sec:construction}
The following theorem is a converse to Theorem~\ref{theorem:squares} which states that if an algebraically periodic surface can be decomposed into squares of various sizes, then its periodic direction field is totally real. 

\begin{thm:1.15-2} \label{totrealK} Let $K$ be a totally real, algebraic extension field of $\QQ$. Then there exists a surface $S$ that arises from a right-angle, square-tileable billiard table so that $S$ is completely algebraically periodic with periodic direction field $K$.  \end{thm:1.15-2}

\begin{proof} Let $K=\QQ(\lambda)$. Since $\lambda$ is totally real, (\cite{Estes}) implies that $\lambda$ is the eigenvalue of an integral, symmetric $n \times n$ matrix $A$. If $u$ is an eigenvector of $A$ corresponding to $\lambda$, then since $A$ is symmetric, $u$ is also an eigenvector for $A^t$ corresponding to $\lambda$. Then if we define
$$
J= \sum_{i=1}^{n} \begin{bmatrix} u_i \\ 0 \end{bmatrix} \wedge \begin{bmatrix} 0 \\ u_i \end{bmatrix}
$$
as in the statement of Theorem~\ref{thm:integral}, the theorem shows that the periodic direction field of $J$ is $K$. Thus, if we can construct a surface $S$ so that $J=J(S)$, $S$ will be have periodic direction field $K$. We begin by constructing the billiard table in the following way. For each term 
$$
\begin{bmatrix} u_i \\ 0 \end{bmatrix} \wedge \begin{bmatrix} 0 \\ u_i \end{bmatrix}
$$
 of $J$, construct a square of side length $u_i$. (Note that we can replace $u_i$ by $-u_i$ without changing the value of $J$. Thus we may assume, without loss of generality, that $u_i >0$.) Glue the squares together along their sides so that the coordinates of the vertices of each square lie in $K$. This may be accomplished, for example, by aligning the squares along the x-axis and gluing them edge to edge. The corresponding translation surface is a four-fold cover of the table, and so $J(S)=4J$ where $J$ is the invariant of the table and $J(S)$ that of the surface. 

It follows from Proposition~\ref{prop:hol criterion} that the holonomy of any segment of $S$ lies in $K^2$.  
By Proposition~\ref{prop:characterization} $S$ is completely algebraically periodic since its periodic direction field is $K$ and its holonomy is contained in $K^2$. 
\end{proof}

Note that the surface $S$ that we have constructed is a connected sum of isogenous tori so this construction provides a converse to Theorem~\ref{thm:isotori}. 

\begin{thm:1.10-2} 
Every number field $K$ is the periodic direction
field for a completely algebraically periodic translation surface that arises from a right-angled billiard table. 
\end{thm:1.10-2}

\begin{proof} 
Let $K$ be a number field. The we can write $K=\QQ(\lambda)$ where
$\lambda>0$ is an algebraic number with minimal polynomial $$p(x)=x^n + a_{n-1}x^{n-1} + \cdots + a_1x + a_0.$$  Let
$$q(x)=\frac{p(x)}{(x-\lambda)}=b_{n-1}x^{n-1} + b_{n-2}x^{n-2} + \cdots + b_1x + b_0.$$
Let $\alpha_j=\lambda^{j}$ and let $\beta_j=\frac{b_{j}}{q'(\lambda)}$. We will construct a surface $S$ so that 
$$J(S)=\sum_i \begin{bmatrix} \alpha_j \\ 0 \end{bmatrix} \wedge \begin{bmatrix} 0 \\ \beta_j \end{bmatrix}.$$
According the Theorem~\ref{thm:J(k)} the periodic direction field of this surface will be $K$.
If $K$ is not totally real, then it may not be the case that $\alpha_j \beta_j > 0$, but it will be the case that $\sum_j \alpha_j \beta_j =1$.

By replacing $\alpha_j$ and $\beta_j$ by $-\alpha_j$ and $-\beta_j$ we may assume that each $\alpha_j >0$. Construct rectangles $R_j\subset\RR^2$ with sides aligned with the axes and with width $\alpha_j$ and height $|\beta_j|$ for $j=1\ldots n$.

Glue together those rectangles $R_j$ for which $\beta_j>0$ so that the vertices of each rectangle are in $K^2$, as in the previous theorem. Now consider those rectangles for which $\beta_j <0$. If one such rectangle $R_j$ has area less than that of a rectangle $R_{j'}$ for which $\beta_{j'}>0$, then cut $R_j$ from $R_{j'}$ in such a way that the vertices are in $K^2$. Otherwise, if $R_j$ has area greater than any rectangle $R_k$ for which $\beta_k >0$, then subdivide $R_j$ into sufficiently small rectangles, keeping their vertices in $K^2$, and remove these smaller rectangles from the positive ones as previously described. 
The result is a planar billiard table whose translation surface has periodic direction field $K$ and holonomy contained in $K^2$. By Proposition~\ref{prop:characterization} the corresponding surface is completely algebraically periodic.
\end{proof} 

Using techniques similar to those in the proof of Theorems~\ref{fieldK}, we now construct a surface $S$ for which the essential holonomy is strictly contained in the absolute holonomy. 

\begin{theorem} 
\label{eh ne ah} There is a translation surface for which the essential holonomy is a proper subspace of the absolute holonomy.
\end{theorem}

\begin{proof}
Let $T$ denote the square billiard table with side length $1$ with sides aligned with the axes.  Note that the direction field, $K$, of this table is $\QQ$. Let $\alpha$ and $\beta$ be irrational, positive real numbers such that $\alpha,\beta < 1$. Let $R$ be a rectangle with sides aligned with the axes and width $\alpha$ and height $\beta$. Remove $R$ from the bottom left corner of $T$ and then glue this removed rectangle to the bottom of the right side of $T$. Call this new table $S$.  
We build a translation surface $S'$ from four copies of $S$. This surface will have genus 3 and $J(S')=4J(S)$.
We have that 

\begin{align*} J(S) & =\begin{bmatrix} 1 \\ 0 \end{bmatrix} \wedge \begin{bmatrix} 0 \\ 1 \end{bmatrix} + \begin{bmatrix} \alpha \\ 0 \end{bmatrix} \wedge \begin{bmatrix} 0 \\ \alpha \end{bmatrix} - \begin{bmatrix} \alpha \\ 0 \end{bmatrix} \wedge \begin{bmatrix} 0 \\ \alpha \end{bmatrix} \\ & = \begin{bmatrix} 1 \\ 0 \end{bmatrix} \wedge \begin{bmatrix} 0 \\ 1 \end{bmatrix} \end{align*} 

Since the periodic direction field of $S$ depends only on $J$ and the value of $J$ is the same as that of the torus the periodic direction field is still $\QQ$.  Since the essential holonomy is the smallest vector space so that $J$ can be written as a sum of wedges of vectors in $EH$, we also know that $EH= \QQ \oplus \QQ$. 
But the absolute holonomy is not contained in $\QQ \oplus \QQ$ since it contains vectors with irrational coordinates. 
\end{proof}

 The surface constructed has genus three but is scissors congruent to the torus. This shows that genus is not a scissors congruence invariant.

Recall that the homological affine group is the subgroup of $\Iso(J)$ which preserves the absolute homology. 

\begin{corollary}
\label{cor:notthesame}
There is a surface $S$ for which $$\Iso(J(S)) \neq \HAG.$$
\end{corollary}

\begin{proof}
Consider the surface constructed above but assume that $\alpha$ and $\beta$ are not rationally related. We have seen that in this case $\Iso(J(S))=SL(2,\QQ)$. On the other hand $\HAG$ is the subgroup of $\Iso(J(S))$ which preserves EH. The intersection of EH with the horizontal axis is $\QQ\oplus\alpha\QQ$ while the intersection of EH with the vertical axis is $\QQ\oplus\beta\QQ$. Since these subspaces of $\RR$ are not isogenous there can be no element of $\HAG$ which takes the x-axis to the y-axis so $\HAG$ is strictly smaller than $SL(2,\QQ)$. 
\end{proof}

\end{section}

\end{document}